\documentclass[10pt]{article}
\usepackage{amssymb,amsmath,amsthm}
\usepackage{CJK}
\usepackage[numbers,sort&compress]{natbib}
\usepackage{xcolor}
\allowdisplaybreaks

\UseRawInputEncoding

\begin{document}
 \begin{CJK}{GBK}{song}
\title{\Large\bf{ Ground state sign-changing homoclinic solutions for a discrete nonlinear $p$-Laplacian equation with logarithmic nonlinearity}}
\date{}
\author {Xin Ou$^1$, \ Xingyong Zhang$^{1,2}$\footnote{Corresponding author, E-mail address: zhangxingyong1@163.com}\\
{\footnotesize $^1$Faculty of Science, Kunming University of Science and Technology,}\\
 {\footnotesize Kunming, Yunnan, 650500, P.R. China.}\\
{\footnotesize $^{2}$Research Center for Mathematics and Interdisciplinary Sciences, Kunming University of Science and Technology,}\\
{\footnotesize Kunming, Yunnan, 650500, P.R. China.}\\
}

 \date{}
 \maketitle

 \begin{center}
 \begin{minipage}{15cm}
 \par
 \small  {\bf Abstract:} By using a direct non-Nehari manifold method from  [X.H. Tang, B.T. Cheng. J. Differ. Equations. 261(2016), 2384-2402.], we obtain an existence result of ground state sign-changing homoclinic solution which only changes sign one times and ground state homoclinic solution for a class of discrete nonlinear $p$-Laplacian equation with logarithmic nonlinearity. Moreover, we prove that the sign-changing ground state energy is larger than twice of the ground state energy.

 \end{minipage}
 \end{center}
  \allowdisplaybreaks
 \vskip2mm
 {\section{Introduction }}
\setcounter{equation}{0}

The existence of solutions for the discrete nonlinear $p$-Laplacian equation by variational methods has always been a hot topic in the last  twenty years and we refer readers to \cite{Mei 2022, Liu 2011, Chen 2013, Zhang 2016, Shi 2016} for example.  Especially, in \cite{Chen 2013}, Chen-Tang considered the following  discrete $p$-Laplacian system:
\begin{eqnarray}\label{ppp1}
   \begin{cases}
   \Delta(\varphi_p(\Delta u(n-1)))-a(n)\varphi_p(u(n))+\nabla W(n,u(n))=0,\;n\in\mathbb{Z},\\
   \lim_{n\rightarrow\pm\infty}u(n)=0,
   \end{cases}
\end{eqnarray}
where $p>1$, $\varphi_p$ is the $p$-Laplace operator, $u\in\mathbb{R}^N$, $a:\mathbb{Z}\to\mathbb{R}$ and $W:\mathbb{Z}\times\mathbb{R}^N\to\mathbb{R}$. When $W(n,x)$ is an odd function in $x$, continuously differentiable and satisfies other suitable conditions,  they obtained that  system has an unbounded sequence of homoclinic solutions using the symmetric mountain pass theorem. When $p=2$, (\ref{ppp1}) reduces to the discrete nonlinear Schr\"{o}dinger(DNLS) equation. The DNLS equation is one of the most important inherently discrete model and plays a crucial role in modeling various phenomena from solid state and condensed matter physics to biology \cite{Flach 2008, Flach 1998, Christodoulides 2003, Hennig 1999}. In recent years, the existence of standing wave solutions for the DNLS equation have attracted some attention (see \cite{Mai 2013, Jia 2017, Yang 2010, Chen 2012, Chen 2018}). Especially, in \cite{Chen 2018}, Chen-Tang-Yu studied the following DNLS equation:
\begin{eqnarray*}
   \begin{cases}
   -\Delta^2u(n-1)+\varepsilon(n)(u(n))-\omega u(n)=f(n,u(n)),n\in\mathbb{Z},\\
   \lim_{n\rightarrow\pm\infty}u(n)=0.
   \end{cases}
\end{eqnarray*}
When $f$ satisfies the superquadratic growth condition and monotonicity condition, using the method in \cite{Cheng 2020} and \cite{Tang 2016}, they obtained that the equation has a ground state solution and a least energy sign-changing solution, which changes sign exactly once. Furthermore, they obtained that  the energy of the sign-changing solution is larger than twice of the ground state energy. Next, we  recall two works \cite{Chang 2023} and \cite{Tang 2016} which inspire our work partially. In \cite{Chang 2023}, Chang-Wang-Yan concerned the following logarithmic Schr\"{o}dinger equation on locally finite graph $G=(V, E)$:
$$
-\Delta u+a(x)u=u\log u^2,\;\; x\in V,
$$
where $a:V\to\mathbb{R}$. When $a$ satisfies is bounded from below and the volume of set $\{x\in V: a(x)\leq M\}$ is finite, they used Nehari manifold method to obtain that the equation has a ground state solution. Moreover, when $a$ is bounded from below  and $1/a(x)$ is Lebesgue integrable function on the set $\{x\in V: a(x)> M_0\}$, they also obtained the equation has a ground state solution by using the mountain pass theorem. In \cite{Tang 2016}, Tang-Cheng investigated the following Kirchhoff type problem:
\begin{eqnarray*}
   \begin{cases}
   -\left(a+b\int_{\Omega}|\nabla u|^2dx\right)\Delta u=f(u),&x\in\Omega,\\
   u=0,&x\in\partial\Omega,
   \end{cases}
\end{eqnarray*}
where $\Omega$ is a bounded domain in $\mathbb{R}^N$, $N=1,2,3$. When $f$ satisfies super-cubic growth and monotonicity condition, they used a new energy inequality, the deformation
lemma, Miranda's theorem and Non-Nehari manifold's method to obtain the existence of ground state solution and sign-changing ground state solution which changes sign exactly once. Furthermore, they obtained that the energy of the sign-changing solution is larger than twice of the ground state energy.
\par
In this paper, inspired by \cite{Chen 2013, Chen 2018, Chang 2023}, we mainly use the method in \cite{Tang 2016} to develop the results in \cite{Chen 2018} to the following discrete nonlinear $p$-Laplacian equation involving logarithmic nonlinearity:
\begin{eqnarray}
\label{eq1}
   \begin{cases}
   -\Delta (a(n-1)\varphi_p(\Delta u(n-1)))+ b(n)\varphi_p(u(n))
=   c(n)|u(n)|^{q-2}u(n)\ln{|u(n)|^r},n\in\mathbb{Z},\\
   \lim_{|n|\rightarrow\infty}u(n)=0,
   \end{cases}
\end{eqnarray}
where $1<p<q$, $\varphi_p(s)=|s|^{p-2}s$ is the $p$-Laplacian operator, $\frac{p}{2}\in\mathbb{N}^*$, $\mathbb{N}^*$ denotes the positive integer set, $a, b, c:\mathbb{Z}\rightarrow (0,+\infty)$, $r\geq 1$, $u:\mathbb Z\to \mathbb R$, and $\Delta u(n)=u(n+1)-u(n)$ is the forward difference operator. Note that the nonlinear term $c(n)|u(n)|^{q-2}u(n)\ln{|u(n)|^r}$  does not satisfy the monotonicity condition in \cite{Chen 2018}. Therefore, the situation we studied is different from that in \cite{Chen 2018} even if $p=2$. There exist two main difficulties in studying the equation (\ref{eq1}). One is that the associated functional $I$ of equation (\ref{eq1}) is not well-defined in $E$, which is caused by the logarithmic nonlinearity,  and the other is that the quasi-linearity  of $p$-Laplacian operator makes establishing energy inequalities difficult and complex. For the first difficulty, we mainly use the idea in \cite{Chang 2023} to establish a well-defined space $\mathcal{D}$, thereby avoiding the case that $\sum_{n\in\mathbb{Z}}c(n)|u(n)|^q\ln{|u(n)|^r}=-\infty$. For the second, we use binomial theorem and combination formula, and then by some careful calculations and analysis, establish some useful energy inequalities. We introduce the following assumptions:
\par
{\it
$(C_1)$ there exists a positive constant $b_0$ such that $b(n)\geq b_0$ for all $n\in \mathbb{Z}$ and $\lim_{|n|\rightarrow +\infty} b(n)=+\infty$;\par
$(C_2)$ there is a positive constant $c_0$ such that $c(n)\leq c_0$ for all $n\in\mathbb{Z}$ and $\sum_{n\in\mathbb{Z}}c(n)<+\infty$.
}
\vskip2mm
Next, we define
$$
V=\Big\{\{u(n)\}_{n\in\mathbb{Z}}:\; u(n)\in\mathbb{R},\;n\in\mathbb{Z}\Big\},
$$
$$
E=\bigg\{u\in V:\sum_{n\in\mathbb{Z}}[a(n)|\Delta u(n)|^p+b(n)|u(n)|^p]<+\infty\bigg\},
$$
and
\begin{eqnarray}
\label{eq6}
\|u\|:=\left(\sum_{n\in\mathbb{Z}}\left[a(n)|\Delta u(n)|^{p}
  +b(n)|u(n)|^{p}\right]\right)^{\frac{1}{p}}.
\end{eqnarray}
Then $E$ is a reflexive Banach space. As usual, let $1<p<+\infty$ and define
$$
l^p(\mathbb{Z},\mathbb{R})=\left\{u\in V:\sum_{n\in\mathbb{Z}}|u(n)|^p<+\infty\right\}
$$
with the norm
\begin{eqnarray*}
\|u\|_{l^p}=\left(\sum_{n\in\mathbb{Z}}|u(n)|^p\right)^\frac{1}{p}.
\end{eqnarray*}
When $p=+\infty$, we define
$$
l^\infty(\mathbb{Z},\mathbb{R})=\left\{u\in V:\sup_{n\in\mathbb{Z}}|u(n)|<+\infty\right\}
$$
with
$$
\|u\|_{l^\infty}=\sup_{n\in\mathbb{Z}}|u(n)|.
$$
Note that equation (\ref{eq1}) is formally related to the energy functional $I: E\rightarrow \mathbb{R}\cup\{+\infty\}$ which is defined by
\begin{eqnarray*}
   I(u)
= \frac{1}{p}\sum_{n\in\mathbb{Z}}\big[a(n)|\Delta u(n)|^p+b(n)|u(n)|^p\big]
  +\frac{r}{q^2}\sum_{n\in\mathbb{Z}}c(n)|u(n)|^q
  -\frac{1}{q}\sum_{n\in\mathbb{Z}}c(n)|u(n)|^{q}\ln {|u(n)|^r}.
\end{eqnarray*}
But the functional $I$ is not well-defined in $E$ (see Appendix 1). We discuss the functional $I$ on the set
\begin{eqnarray*}
  \mathcal{D}
= \Big\{u\in E:\sum_{n\in\mathbb{Z}}c(n)|u(n)|^q\ln{|u(n)|^r}<+\infty\Big\},
\end{eqnarray*}
that is,
\begin{eqnarray}
\label{eq7}
  I(u)
= \frac{1}{p}\|u\|^p+\frac{r}{q^2}\sum_{n\in\mathbb{Z}}c(n)|u(n)|^q
  -\frac{1}{q}\sum_{n\in\mathbb{Z}}c(n)|u(n)|^{q}\ln {|u(n)|^r},\;\;\forall u\in\mathcal{D}.
\end{eqnarray}
Note that
\begin{eqnarray}
\label{eq8}
\lim_{t\rightarrow 0}\frac{t^{q-1}\ln |t|^r}{t^{p-1}}=0&\text {and}&\lim_{t\rightarrow \infty}\frac{t^{q-1}\ln |t|^r}{t^{\zeta-1}}=0
\end{eqnarray}
for all $n\in\mathbb{Z}$, where $\zeta\in(q,+\infty)$. Then by $(C_2)$, for any given $\varepsilon >0$, there exists a positive constant $C_\varepsilon$ such that
\begin{eqnarray}
\label{eq9}
c(n)|t|^{q-1}\big|\ln |t|^r\big|\leq c(n)\varepsilon|t|^{p-1}+c(n)C_\varepsilon|t|^{\zeta-1}\leq c_0\varepsilon|t|^{p-1}+c_0C_\varepsilon|t|^{\zeta-1},& \forall t\in\mathbb{R},\ \forall n\in\mathbb{Z}.
\end{eqnarray}
Then $\mathcal{D}$ is the closed subspace of $E$,  $I\in C^1(\mathcal{D},\mathbb{R})$ and
\begin{eqnarray}
\label{eq10}
      \langle I'(u),v \rangle
& =&  \sum_{n\in\mathbb{Z}}\left[a(n)|\Delta u(n)|^{p-2}\Delta u(n)\Delta v(n)+b(n)|u(n)|^{p-2}u(n)v(n)\right]\nonumber\\
&  &   -\sum_{n\in\mathbb{Z}}c(n)|u(n)|^{q-2}u(n)v(n)\ln {|u(n)|^r},\;\;\forall u,v\in\mathcal{D}.
\end{eqnarray}
Using Abel's partial summation formula (also known as Abel's transformation) in \cite{Abel} and the definition of $\Delta u(n)$, we have
\begin{eqnarray*}
\label{eq11}
  \sum_{n\in\mathbb{Z}}a(n)|\Delta u(n)|^{p-2}\Delta u(n)\Delta v(n)
= -\sum_{n\in\mathbb{Z}}\Delta\Big(a(n-1)\varphi_p\big(\Delta u(n-1)\big)\Big)v(n),
\end{eqnarray*}
which implies that
\begin{eqnarray*}
      \langle I'(u),v \rangle
& =&  \sum_{n\in\mathbb{Z}}\left[-\Delta\Big(a(n-1)\varphi_p\big(\Delta u(n-1)\big)\Big)v(n)
       +b(n)|u(n)|^{p-2}u(n)v(n)\right]\nonumber\\
&  &   -\sum_{n\in\mathbb{Z}}c(n)|u(n)|^{q-2}u(n)v(n)\ln {|u(n)|^r},\;\;\forall u,v\in\mathcal{D}.
\end{eqnarray*}
According to the above equations, we can derive that $\langle I'(u),v\rangle=0$ for any $v\in \mathcal{D}$ if and only if
$$
  -\Delta\Big(a(n-1)\varphi_p\big(\Delta u(n-1)\big)\Big)+b(n)|u(n)|^{p-2}u(n)
=  c(n)|u(n)|^{q-2}u(n)\ln {|u(n)|^r}.
$$
Therefore, it is easy to see that the critical points of $I$ in $\mathcal{D}$ are solutions of equation (\ref{eq1}). Furthermore, if $u\in\mathcal{D}$ is a solution of equation (\ref{eq1}) and $u^\pm\neq 0$, then $u$ is a sign-changing solution of equation (\ref{eq1}), where
$$
u^+(n):=\max\{u(n),0\}\;\;\text{and}\;\;u^-(n):=\min\{u(n),0\}.
$$
\par
To be precise, we obtain the following results.
\par

\vskip2mm
\noindent
{\bf Theorem 1.1.} {\it Assume that $(C_1)$ and $(C_2)$ hold. Then problem} (\ref{eq1}) {\it has a sign-changing solution $u_0\in \mathcal{M}$ such that $I(u_0)=\inf_{\mathcal{M}}I:=m_*>0$ and $u_0$ only changes the sign one times, where}
\begin{eqnarray*}
  \mathcal{M}
= \big\{u\in \mathcal{D}:u^{\pm}\neq 0,\;\;\langle I'(u),u^{+}\rangle=0\;\;and\;\;\langle I'(u),u^{-}\rangle=0\big\}.
\end{eqnarray*}
\vskip2mm
\noindent
{\bf Theorem 1.2.} {\it Assume that $(C_1)$ and $(C_2)$ hold. Then problem} (\ref{eq1}) {\it has a solution $\bar{u}\in \mathcal{N}$ such that $I(\bar{u})=\inf_{\mathcal{N}}I:=c_*>0$, where
\begin{eqnarray*}
  \mathcal{N}
= \big\{u\in \mathcal{D}:u\neq 0\;and\; \langle I'(u),u\rangle=0\big\}.
\end{eqnarray*}
In addition, $m_*\geq 2c_*$.}

\vskip2mm
{\section{Preliminaries}}
\setcounter{equation}{0}

In this section, we provide some lemmas which play some important roles in the proofs of our results.
\vskip2mm
\noindent
{\bf Lemma 2.1.} {\it Assume that $(C_1)$ holds. Then $\mathcal{D}$ is continuously embedded into $l^\kappa(\mathbb{Z},\mathbb{R})$ for any $p\leq\kappa\leq+\infty$, that is, for all $u\in\mathcal{D}$,
\begin{eqnarray}
\label{Eq1}
\|u\|_{l^\kappa}\leq b_0^{-\frac{1}{p}}\|u\|,
\end{eqnarray}
Moreover, $\mathcal{D}$ is compactly embedded in $l^\kappa(\mathbb{Z},\mathbb{R})$ for any $p\leq\kappa\leq+\infty$.}\\
{\bf Proof.} For any $u\in\mathcal{D}$,when $\kappa=p$, there holds
$$
     \|u\|
\geq \left(\sum_{n\in\mathbb{Z}}b(n)|u(n)|^p\right)^{\frac{1}{p}}
\geq b_0^{\frac{1}{p}}\left(\sum_{n\in\mathbb{Z}}|u(n)|^p\right)^{\frac{1}{p}}
 =    b_0^{\frac{1}{p}}\|u\|_{l^p}.
$$
When $\kappa=+\infty$, we can also obtain that
\begin{eqnarray}\label{k1}
     \|u\|
\geq \left(\sum_{n\in\mathbb{Z}}b(n)|u(n)|^p\right)^{\frac{1}{p}}
\geq b_0^{\frac{1}{p}}\left(\sum_{n\in\mathbb{Z}}|u(n)|^p\right)^{\frac{1}{p}}
\geq b_0^{\frac{1}{p}}\left(\sup_{n\in\mathbb{Z}}|u(n)|^p\right)^{\frac{1}{p}}
 =   b_0^{\frac{1}{p}}\|u\|_{l^\infty},
\end{eqnarray}
For any $p<\kappa<+\infty$, it follows from (\ref{k1}) that
\begin{eqnarray*}
&    & \|u\|_{l^\kappa}^\kappa=\sum_{n\in\mathbb{Z}}|u(n)|^\kappa=\sum_{n\in\mathbb{Z}}|u(n)|^{\kappa-p}|u(n)|^p\\
&\leq& \|u\|^{\kappa-p}_{l^\infty}\sum_{n\in\mathbb{Z}}|u(n)|^p
               =\|u\|^{\kappa-p}_{l^\infty}\sum_{n\in\mathbb{Z}}\frac{1}{b(n)}b(n)|u(n)|^p\\
&\leq&  b_0^{-\frac{\kappa-p}{p}}\|u\|^{\kappa-p}\frac{1}{b_0}\|u\|^p=b_0^{-\frac{\kappa}{p}}\|u\|^{\kappa}.
\end{eqnarray*}
So, (\ref{Eq1}) holds.
\par
Next, we prove that the embeddings are also compact. Suppose that $\{u_k\}$ is a bounded sequence in $\mathcal{D}$. Then there is a subsequence of $\{u_k\}$, still denoted by $\{u_k\}$, such that $u_k\rightharpoonup u$ weakly in $\mathcal{D}$ for some point $u\in\mathcal{D}$. In particular,
$$
\lim_{k \rightarrow +\infty}\sum_{n\in\mathbb{Z}}u_k\varphi=\sum_{n\in\mathbb{Z}}u\varphi,
$$
where  $\varphi\in \mathcal{D}$ is defined by
\begin{eqnarray*}
     \varphi(m)=
  \begin{cases}
   1,\;\;m=n,\\
   0,\;\;m\neq n
  \end{cases}
\end{eqnarray*}
for any fixed $n$. Thus, we have
\begin{eqnarray}
\label{Equ1}
\lim_{k\rightarrow +\infty} u_k(n)=u(n)\;\;\text{for any fixed}\;\; n\in\mathbb{Z}.
\end{eqnarray}
We now prove $u_k\rightarrow u$ in $l^\kappa(\mathbb{Z},\mathbb{R})$ for all $p\leq\kappa\leq +\infty$. When $\kappa=p$, since $u\in\mathcal{D}$, according to the boundedness of $\{u_k\}$ and the definition of $\|\cdot\|$, there appears a positive constant $\delta_0$ such that
\begin{eqnarray*}\label{k2}
\sum_{n\in\mathbb{Z}}b(n)|u_k(n)-u(n)|^p\leq\delta_0.
\end{eqnarray*}
For any given positive constant $\varepsilon_1$, there is a $n_0\in\mathbb Z$ such that $\frac{1}{b(n)}<\varepsilon_1\;\;\text{as }|n|>n_0$.
Therefore, we can obtain that
\begin{eqnarray}
\label{Equ2}
\sum_{|n|>n_0}|u_k(n)-u(n)|^p=\sum_{|n|>n_0}\frac{1}{b(n)}b(n)|u_k(n)-u(n)|^p\leq\varepsilon_1\delta_0.
\end{eqnarray}
On the other hand, $(\ref{Equ1})$ implies that $\lim\limits_{k\rightarrow +\infty}\sum\limits_{|n|\le n_0}|u_k(n)-u(n)|^p=0$ since $\{n\in\mathbb{Z}:|n|\le n_0\}$ is a finite set. Then, according to the arbitrariness of $\varepsilon_1$ and $(\ref{Equ2})$, we have
\begin{small}
\begin{eqnarray}
\label{Equ3}
 \lim_{k\rightarrow +\infty}\sum_{n\in\mathbb{Z}}|u_k(n)-u(n)|^p=0.
\end{eqnarray}
\end{small}
For $\kappa=+\infty$, according to the definition of $\|\cdot\|_{l^\infty}$ and $(\ref{Equ3})$, as $k\rightarrow +\infty$, we have
\begin{eqnarray}
\label{Equ4}
       \|u_k-u\|_{l^\infty}^p
\leq \sum_{n\in\mathbb{Z}}|u_k(n)-u(n)|^p\rightarrow 0,
\end{eqnarray}
and for $p<\kappa<+\infty$, by (\ref{Equ3}) and (\ref{Equ4}), there exists
\begin{eqnarray}
\label{Equ5}
&   & \|u_k-u\|_{l^\kappa}^\kappa=\sum_{n\in\mathbb{Z}}|u_k(n)-u(n)|^{\kappa-p}|u_k(n)-u(n)|^p\nonumber\\
&\le& \|u_k-u\|_{l^\infty}^{\kappa-p}\sum_{n\in\mathbb{Z}}|u_k(n)-u(n)|^p
                  \rightarrow 0,\;\;\text{as}\;\;k\rightarrow +\infty.
\end{eqnarray}
Consequently, by (\ref{Equ3}), (\ref{Equ4}) and (\ref{Equ5}), we can derive that $u_k\rightarrow u$ in $l^\kappa(\mathbb{Z},\mathbb{R})$ for all $p\leq\kappa\leq +\infty$.
\qed
\vskip2mm
\noindent
{\bf Proposition 2.1.} {\it For all $\frac{p}{2}\in\mathbb{N}^*$ and $u\in\mathcal{D}$, there hold
\begin{small}
\begin{eqnarray*}
      I(u)
& = & I(u^+)+I(u^-)+\frac{1}{p}\sum_{n\in\mathbb{Z}}a(n)\Big|\sum_{i=1}^{\frac{p}{2}-1}\sum_{j=0}^iC_{\frac{p}{2}}^iC_i^j2^{i-j}
           (\Delta u^+(n))^{p-(i+j)}(\Delta u^-(n))^{i+j}\Big|\nonumber\\
&   & +\frac{1}{p}\sum_{n\in\mathbb{Z}}a(n)\Big|\sum_{j=0}^{\frac{p}{2}-1}C_{\frac{p}{2}}^j2^{\frac{p}{2}-j}
           (\Delta u^+(n))^{\frac{p}{2}-j}(\Delta u^-(n))^{\frac{p}{2}+j}\Big|,
\end{eqnarray*}
\begin{eqnarray*}
     \langle I'(u),u^+\rangle
&=&  \langle I'(u^+),u^+\rangle
       +\sum_{n\in\mathbb{Z}}a(n)\Big|\sum_{i=1}^{\frac{p}{2}-1}\sum_{j=0}^i C_{\frac{p}{2}-1}^i C_i^j2^{i-j}
              (\Delta u^+(n))^{p-(i+j)}(\Delta u^-(n))^{i+j}\Big|\nonumber\\
&  & +\sum_{n\in\mathbb{Z}}a(n)\Big|\sum_{i=1}^{\frac{p}{2}}\sum_{j=0}^{i-1} C_{\frac{p}{2}-1}^{i-1} C_{i-1}^{j}2^{i-1-j}
              (\Delta u^+(n))^{p-(i+j)}(\Delta u^-(n))^{i+j}\Big|,
\end{eqnarray*}
and
\begin{eqnarray*}
      \langle I'(u),u^-\rangle
&=&   \langle I'(u^-),u^-\rangle
       +\sum_{n\in\mathbb{Z}}a(n)\Big|\sum_{i=1}^{\frac{p}{2}-1}\sum_{j=1}^{i} C_{\frac{p}{2}-1}^{i-1} C_{i-1}^{j-1}2^{i-j}
               (\Delta u^+(n))^{p-(i+j)}(\Delta u^-(n))^{i+j}\Big|\nonumber\\
&  &  +\sum_{n\in\mathbb{Z}}a(n)\Big|\sum_{j=1}^{\frac{p}{2}-1}C_{\frac{p}{2}-1}^{j-1}2^{\frac{p}{2}-j}
               (\Delta u^+(n))^{\frac{p}{2}-j}(\Delta u^-(n))^{\frac{p}{2}+j}\Big|\nonumber\\
&  &  +\sum_{n\in\mathbb{Z}}a(n)\Big|\sum_{i=1}^{\frac{p}{2}}\sum_{j=0}^{i-1} C_{\frac{p}{2}-1}^{i-1} C_{i-1}^{j}2^{i-1-j}
              (\Delta u^+(n))^{p-(i+j)}(\Delta u^-(n))^{i+j}\Big|.
\end{eqnarray*}
\end{small}
}
{\bf Proof.} Let
\begin{eqnarray*}
\mathbb{Z}_1:=\{n\in\mathbb{Z}:u(n)\geq0\}\;\;\text{and}\;\;\mathbb{Z}_2:=\{n\in\mathbb{Z}:u(n)<0\}.
\end{eqnarray*}
Note that $\Delta u^+(n)\Delta u^-(n)=-u^+(n+1)u^-(n)-u^+(n)u^-(n+1)\geq0$. Then according to the definition of $\|\cdot\|$, Appendix 1 below and binomial theorem, we have
\begin{small}
\begin{eqnarray}
\label{Eq2}
       \|u\|^p
& = &  \sum_{n\in\mathbb{Z}}\big[a(n)|\Delta u^+(n)+\Delta u^-(n)|^p+b(n)|u^+(n)+u^-(n)|^p\big]\nonumber\\
& = &  \sum_{n\in\mathbb{Z}}\left[a(n)\Big|(\Delta u^+(n))^2+2\Delta u^+(n)\Delta u^-(n)+(\Delta u^-(n))^2\Big|^{\frac{p}{2}}
         +b(n)|u^+(n)+u^-(n)|^p\right]\nonumber\\
& = &  \sum_{n\in\mathbb{Z}}a(n)\Big|\sum_{i=0}^{\frac{p}{2}}C_{\frac{p}{2}}^i(\Delta u^+(n))^{p-2i}
          \sum_{j=0}^iC_i^j(2\Delta u^+(n)\Delta u^-(n))^{i-j}(\Delta u^-(n))^{2j}\Big|  +\sum_{n\in\mathbb{Z}}b(n)|u^+(n)+u^-(n)|^p\nonumber\\
& = &  \sum_{n\in\mathbb{Z}}a(n)\Big|\sum_{i=0}^{\frac{p}{2}}\sum_{j=0}^iC_{\frac{p}{2}}^iC_i^j2^{i-j}(\Delta u^+(n))^{p-(i+j)}
          (\Delta u^-(n))^{i+j}\Big|+\sum_{n\in\mathbb{Z}_1}b(n)|u^+(n)|^p+\sum_{n\in\mathbb{Z}_2}b(n)|u^-(n)|^p\nonumber\\
& = &  \sum_{n\in\mathbb{Z}}[a(n)|\Delta u^+(n)|^p+b(n)|u^+(n)|^p]
          +\sum_{n\in\mathbb{Z}}a(n)\Big|\sum_{i=1}^{\frac{p}{2}-1}\sum_{j=0}^iC_{\frac{p}{2}}^iC_i^j2^{i-j}
           (\Delta u^+(n))^{p-(i+j)}(\Delta u^-(n))^{i+j}\Big|\nonumber\\
&¡¡ &  +\sum_{n\in\mathbb{Z}}[a(n)|\Delta u^-(n)|^p+b(n)|u^-(n)|^p]
          +\sum_{n\in\mathbb{Z}}a(n)\Big|\sum_{j=0}^{\frac{p}{2}-1}C_{\frac{p}{2}}^j2^{\frac{p}{2}-j}
           (\Delta u^+(n))^{\frac{p}{2}-j}(\Delta u^-(n))^{\frac{p}{2}+j}\Big|\nonumber\\
& = &  \|u^+\|^p+\|u^-\|^p+\sum_{n\in\mathbb{Z}}a(n)\Big|\sum_{i=1}^{\frac{p}{2}-1}\sum_{j=0}^iC_{\frac{p}{2}}^iC_i^j2^{i-j}
           (\Delta u^+(n))^{p-(i+j)}(\Delta u^-(n))^{i+j}\Big|\nonumber\\
&   &  +\sum_{n\in\mathbb{Z}}a(n)\Big|\sum_{j=0}^{\frac{p}{2}-1}C_{\frac{p}{2}}^j2^{\frac{p}{2}-j}
           (\Delta u^+(n))^{\frac{p}{2}-j}(\Delta u^-(n))^{\frac{p}{2}+j}\Big|.
\end{eqnarray}
\end{small}
Similarly, we have
\begin{eqnarray}
\label{Eq3}
&   &  \sum_{n\in\mathbb{Z}}[a(n)|\Delta u(n)|^{p-2}\Delta u(n)\Delta u^+(n)+b(n)|u(n)|^{p-2}u(n)u^+(n)]\nonumber\\
& = &  \sum_{n\in\mathbb{Z}}\left[a(n)|\Delta u^+(n)+\Delta u^-(n)|^{p-2}(\Delta u^+(n)+\Delta u^-(n))\Delta u^+(n)
         +b(n)|u(n)|^{p-2}u(n)u^+(n)\right]\nonumber\\
& = &  \sum_{n\in\mathbb{Z}}\left[a(n)|(\Delta u^+(n))^2+2\Delta u^+(n)\Delta u^-(n)+(\Delta u^-(n))^2|^{\frac{p-2}{2}}
                 (\Delta u^+(n))^2\right]+\sum_{n\in\mathbb{Z}_1}b(n)|u^+(n)|^{p}\nonumber\\
&   &  +\sum_{n\in\mathbb{Z}}\left[a(n)|(\Delta u^+(n))^2+2\Delta u^+(n)\Delta u^-(n)+(\Delta u^-(n))^2|^{\frac{p-2}{2}}
                 \Delta u^-(n)\Delta u^+(n)\right]\nonumber\\
& = &  \sum_{n\in\mathbb{Z}}\left[a(n)\Big|\sum_{i=0}^{\frac{p}{2}-1}\sum_{j=0}^i C_{\frac{p}{2}-1}^i C_i^j2^{i-j}
(\Delta u^+(n))^{p-2-(i+j)}(\Delta u^-(n))^{i+j}\Big|(\Delta u^+(n))^2\right]+\sum_{n\in\mathbb{Z}_1}b(n)|u^+(n)|^{p}\nonumber\\
&  &   +\sum_{n\in\mathbb{Z}}\left[a(n)\Big|\sum_{i=0}^{\frac{p}{2}-1}\sum_{j=0}^i C_{\frac{p}{2}-1}^i C_i^j2^{i-j}
                 (\Delta u^+(n))^{p-2-(i+j)}(\Delta u^-(n))^{i+j}\Big|\Delta u^-(n)\Delta u^+(n)\right]\nonumber\\
& = &  \|u^+\|^p+\sum_{n\in\mathbb{Z}}a(n)\Big|\sum_{i=1}^{\frac{p}{2}-1}\sum_{j=0}^i C_{\frac{p}{2}-1}^i C_i^j2^{i-j}
                 (\Delta u^+(n))^{p-(i+j)}(\Delta u^-(n))^{i+j}\Big|\nonumber\\
&  &   +\sum_{n\in\mathbb{Z}}a(n)\Big|\sum_{i=1}^{\frac{p}{2}}\sum_{j=0}^{i-1} C_{\frac{p}{2}-1}^{i-1} C_{i-1}^{j}2^{i-1-j}
              (\Delta u^+(n))^{p-(i+j)}(\Delta u^-(n))^{i+j}\Big|,
\end{eqnarray}
and
\begin{eqnarray}
\label{Eq4}
&   &  \sum_{n\in\mathbb{Z}}[a(n)|\Delta u(n)|^{p-2}\Delta u(n)\Delta u^-(n)+b(n)|u(n)|^{p-2}u(n)u^-(n)]\nonumber\\
& = &  \sum_{n\in\mathbb{Z}}\left[a(n)\Big|\sum_{i=0}^{\frac{p}{2}-1}\sum_{j=0}^i C_{\frac{p}{2}-1}^i C_i^j2^{i-j}
(\Delta u^+(n))^{p-2-(i+j)}(\Delta u^-(n))^{i+j}\Big|(\Delta u^-(n))^2\right]+\sum_{n\in\mathbb{Z}_2}b(n)|u^-(n)|^{p}\nonumber\\
&  &   +\sum_{n\in\mathbb{Z}}\left[a(n)\Big|\sum_{i=0}^{\frac{p}{2}-1}\sum_{j=0}^i C_{\frac{p}{2}-1}^i C_i^j2^{i-j}
                       (\Delta u^+(n))^{p-2-(i+j)}(\Delta u^-(n))^{i+j}\Big|\Delta u^+(n)\Delta u^-(n)\right]\nonumber\\
& = &  \|u^-\|^p+\sum_{n\in\mathbb{Z}}a(n)\Big|\sum_{i=0}^{\frac{p}{2}-2}\sum_{j=0}^i C_{\frac{p}{2}-1}^i C_i^j2^{i-j}
                       (\Delta u^+(n))^{p-(i+j+2)}(\Delta u^-(n))^{i+j+2}\Big|\nonumber\\
&  &   +\sum_{n\in\mathbb{Z}}a(n)\Big|\sum_{j=0}^{\frac{p}{2}-2}C_{\frac{p}{2}-1}^j2^{{\frac{p}{2}-1}-j}
                       (\Delta u^+(n))^{p-({\frac{p}{2}-1}+j+2)}(\Delta u^-(n))^{{\frac{p}{2}-1}+j+2}\Big|\nonumber\\
&  &   +\sum_{n\in\mathbb{Z}}a(n)\Big|\sum_{i=0}^{\frac{p}{2}-1}\sum_{j=0}^i C_{\frac{p}{2}-1}^i C_i^j2^{i-j}
                       (\Delta u^+(n))^{p-1-(i+j)}(\Delta u^-(n))^{i+j+1}\Big|\nonumber\\
& = &  \|u^-\|^p+\sum_{n\in\mathbb{Z}}a(n)\Big|\sum_{i=1}^{\frac{p}{2}-1}\sum_{j=1}^{i} C_{\frac{p}{2}-1}^{i-1}
                       C_{i-1}^{j-1}2^{i-j}(\Delta u^+(n))^{p-(i+j)}(\Delta u^-(n))^{i+j}\Big|\nonumber\\
&  &   +\sum_{n\in\mathbb{Z}}a(n)\Big|\sum_{j=1}^{\frac{p}{2}-1}C_{\frac{p}{2}-1}^{j-1}2^{\frac{p}{2}-j}
                       (\Delta u^+(n))^{\frac{p}{2}-j}(\Delta u^-(n))^{\frac{p}{2}+j}\Big|\nonumber\\
&  &   +\sum_{n\in\mathbb{Z}}a(n)\Big|\sum_{i=1}^{\frac{p}{2}}\sum_{j=0}^{i-1} C_{\frac{p}{2}-1}^{i-1} C_{i-1}^{j}2^{i-1-j}
              (\Delta u^+(n))^{p-(i+j)}(\Delta u^-(n))^{i+j}\Big|.
\end{eqnarray}
By (\ref{eq7}), (\ref{eq10}), (\ref{Eq2}), (\ref{Eq3}) and (\ref{Eq4}), it is easy to see that the conclusions hold.
\qed
\par
Next, we establish an inequality associated to $I(u)$, $I(su^++tu^-)$, $\langle I'(u),u^+\rangle$ and $\langle I'(u),u^-\rangle$.
\vskip2mm
\noindent
{\bf Lemma 2.2.} {\it Assume that $(C_1)$ and $(C_2)$ hold. For all $u\in\mathcal{D}$ and $s,t\ge 0$, there exists
{\small \begin{eqnarray}
\label{Eq5}
        I(u)
&\geq& I(su^++tu^-)+\frac{1-s^{q}}{q}\langle I'(u),u^+\rangle
           +\frac{1-t^{q}}{q}\langle I'(u),u^-\rangle+\left(\frac{1-s^p}{p}-\frac{1-s^{q}}{q}\right)\|u^+\|^p\nonumber\\
&    & +\left(\frac{1-t^p}{p}-\frac{1-t^{q}}{q}\right)\|u^-\|^p
           +\sum_{n\in\mathbb{Z}}a(n)\Big|\sum_{i=1}^{\frac{p}{2}-1}\sum_{j=1}^{i-1}2^{i-j}(\Delta u^+(n))^{p-(i+j)}
                   (\Delta u^-(n))^{i+j}\Big|\Theta,
\end{eqnarray}}
where $\Theta=\frac{2s^pC_{\frac{p}{2}-1}^i C_i^j+s^pC_{\frac{p}{2}-1}^{i-1}C_{i-1}^j+2t^pC_{\frac{p}{2}-1}^{i-1}C_{i-1}^{j-1}
+t^pC_{\frac{p}{2}-1}^{i-1}C_{i-1}^j-2s^{p-(i+j)}t^{i+j}C_{\frac{p}{2}}^i C_i^j}{2p}\geq0$.}\\
{\bf Proof.} It is easy to see that (\ref{Eq5}) holds for $u=0$. Next, we let $u\neq 0$. According to Appendix 1 in \cite{Zhang 2023}, there holds
\begin{eqnarray}
\label{Eq6}
r(1-\tau^{q})+q\tau^{q}\ln\tau^r> 0,\;\;\forall \tau\in(0,1)\cup(1,+\infty).
\end{eqnarray}
For $u\in \mathcal{D}\backslash\{0\}$ and all $s,t\ge 0$, we have
\begin{eqnarray}
\label{Eq7}
&    & \sum_{n\in\mathbb{Z}}|su^+(n)+tu^-(n)|^{q}\ln |su^+(n)+tu^-(n)|^r\nonumber\\
& =  & \sum_{n\in\mathbb{Z}_1}|su^+(n)+tu^-(n)|^{q}\ln |su^+(n)+tu^-(n)|^r\nonumber\\
&    & +\sum_{n\in\mathbb{Z}_2}|su^+(n)+tu^-(n)|^{q}\ln |su^+(n)+tu^-(n)|^r\nonumber\\
& =  & \sum_{n\in\mathbb{Z}_1}|su^+(n)|^{q}\ln |su^+(n)|^r
        +\sum_{n\in\mathbb{Z}_2}|tu^-(n)|^{q}\ln |tu^-(n)|^r\nonumber\\
& =  & \sum_{n\in\mathbb{Z}}\big[|su^+(n)|^{q}(\ln |u^+(n)|^r+\ln s^r)+|tu^-(n)|^{q}(\ln |u^-(n)|^r+\ln t^r)\big].
\end{eqnarray}
In virtue of Appendix 2 in \cite{Zhang 2023},  the function $f(x)=\frac{1-a^x}{x}$ is strictly monotonically decreasing on $(0,+\infty)$ for $a>0$ and $a\neq 1$. Then by (\ref{eq7}), (\ref{Eq2}), (\ref{Eq7}), (\ref{Eq3}), (\ref{Eq4}), (\ref{eq10}), (\ref{Eq6}) and Appendix 2 below, we can derive the following inequality
\begin{eqnarray}
\label{Eq8}
&    & I(u)-I(su^++tu^-)\nonumber\\
& =  & \frac{1}{p}\left(\|u\|^p-\|su^++tu^-\|^p\right)
          +\frac{r}{q^2}\sum_{n\in\mathbb{Z}}c(n)\left(|u(n)|^q-|su^+(n)+ tu^-(n)|^q\right)\nonumber\\
&    & -\frac{1}{q}\sum_{n\in\mathbb{Z}}c(n)
          \left(|u(n)|^q\ln{|u(n)|^r-|su^+(n)+tu^-(n)|^q\ln{|su^+(n)+tu^-(n)|^r}}\right)\nonumber\\
& =  & \frac{1-s^p}{p}\|u^+\|^p+\frac{1-t^p}{p}\|u^-\|^p
          +\frac{r(1-s^q)}{q^2}\sum_{n\in\mathbb{Z}}c(n)|u^+(n)|^q+\frac{r(1-t^q)}{q^2}\sum_{n\in\mathbb{Z}}c(n)|u^-(n)|^q\nonumber\\
&    & -\frac{1-s^q}{q}\sum_{n\in\mathbb{Z}}c(n)|u^+(n)|^q\ln{|u^+(n)|^r}
          -\frac{1-t^q}{q}\sum_{n\in\mathbb{Z}}c(n)|u^-(n)|^q\ln{|u^-(n)|^r}\nonumber\\
&    & +\frac{1}{q}\sum_{n\in\mathbb{Z}}c(n)|su^+(n)|^q\ln{s^r}
          +\frac{1}{q}\sum_{n\in\mathbb{Z}}c(n)|tu^-(n)|^q\ln{t^r}\nonumber\\
&    & +\frac{1}{p}\sum_{n\in\mathbb{Z}}a(n)\Big|\sum_{i=1}^{\frac{p}{2}-1}\sum_{j=0}^iC_{\frac{p}{2}}^iC_i^j2^{i-j}
           (\Delta u^+(n))^{p-(i+j)}(\Delta u^-(n))^{i+j}\Big|\nonumber\\
&   &  +\frac{1}{p}\sum_{n\in\mathbb{Z}}a(n)\Big|\sum_{j=0}^{\frac{p}{2}-1}C_{\frac{p}{2}}^j2^{\frac{p}{2}-j}
           (\Delta u^+(n))^{\frac{p}{2}-j}(\Delta u^-(n))^{\frac{p}{2}+j}\Big|\nonumber\\
&   &  -\frac{1}{p}\sum_{n\in\mathbb{Z}}a(n)\Big|\sum_{i=1}^{\frac{p}{2}-1}\sum_{j=0}^iC_{\frac{p}{2}}^iC_i^j2^{i-j}
           (\Delta su^+(n))^{p-(i+j)}(\Delta tu^-(n))^{i+j}\Big|\nonumber\\
&   &  -\frac{1}{p}\sum_{n\in\mathbb{Z}}a(n)\Big|\sum_{j=0}^{\frac{p}{2}-1}C_{\frac{p}{2}}^j2^{\frac{p}{2}-j}
           (\Delta su^+(n))^{\frac{p}{2}-j}(\Delta tu^-(n))^{\frac{p}{2}+j}\Big|\nonumber\\
& =  & \frac{1-s^q}{q}\left[\langle I'(u),u^+\rangle-\sum_{n\in\mathbb{Z}}a(n)|\Delta u(n)|^{p-2}\Delta u(n)\Delta u^+(n)
          -\sum_{n\in\mathbb{Z}}b(n)|u^+(n)|^p\right]\nonumber\\
&    & +\frac{1-t^q}{q}\left[\langle I'(u),u^-\rangle-\sum_{n\in\mathbb{Z}}a(n)|\Delta u(n)|^{p-2}\Delta u(n)\Delta u^-(n)
          -\sum_{n\in\mathbb{Z}}b(n)|u^-(n)|^p\right]\nonumber\\
&    & +\frac{1-s^p}{p}\|u^+\|^p+\frac{1-t^p}{p}\|u^-\|^p
          +\frac{r(1-s^q)}{q^2}\sum_{n\in\mathbb{Z}}c(n)|u^+(n)|^q+\frac{r(1-t^q)}{q^2}\sum_{n\in\mathbb{Z}}c(n)|u^-(n)|^q\nonumber\\
&    & +\frac{1}{q}\sum_{n\in\mathbb{Z}}c(n)|su^+(n)|^q\ln{s^r}
          +\frac{1}{q}\sum_{n\in\mathbb{Z}}c(n)|tu^-(n)|^q\ln{t^r}\nonumber\\
&    & +\frac{1}{p}\sum_{n\in\mathbb{Z}}a(n)\Big|\sum_{i=1}^{\frac{p}{2}-1}\sum_{j=0}^iC_{\frac{p}{2}}^iC_i^j2^{i-j}
           (\Delta u^+(n))^{p-(i+j)}(\Delta u^-(n))^{i+j}\Big|\nonumber\\
&   &  +\frac{1}{p}\sum_{n\in\mathbb{Z}}a(n)\Big|\sum_{j=0}^{\frac{p}{2}-1}C_{\frac{p}{2}}^j2^{\frac{p}{2}-j}
           (\Delta u^+(n))^{\frac{p}{2}-j}(\Delta u^-(n))^{\frac{p}{2}+j}\Big|\nonumber\\
&   &  -\frac{1}{p}\sum_{n\in\mathbb{Z}}a(n)\Big|\sum_{i=1}^{\frac{p}{2}-1}\sum_{j=0}^iC_{\frac{p}{2}}^iC_i^j2^{i-j}
           (\Delta su^+(n))^{p-(i+j)}(\Delta tu^-(n))^{i+j}\Big|\nonumber\\
&   &  -\frac{1}{p}\sum_{n\in\mathbb{Z}}a(n)\Big|\sum_{j=0}^{\frac{p}{2}-1}C_{\frac{p}{2}}^j2^{\frac{p}{2}-j}
           (\Delta su^+(n))^{\frac{p}{2}-j}(\Delta tu^-(n))^{\frac{p}{2}+j}\Big|\nonumber\\
& =  & \frac{1-s^{q}}{q}\langle I'(u),u^+\rangle+\frac{1-t^{q}}{q}\langle I'(u),u^-\rangle
       +\left(\frac{1-s^p}{p}-\frac{1-s^{q}}{q}\right)\|u^+\|^p+\left(\frac{1-t^p}{p}-\frac{1-t^{q}}{q}\right)\|u^-\|^p\nonumber\\
&    & +\frac{r(1-s^{q})+qs^{q}\ln s^r}{q^2}\sum_{n\in\mathbb{Z}}c(n)|u^+(n)|^q
         +\frac{r(1-t^{q})+qt^{q}\ln t^r}{q^2}\sum_{n\in\mathbb{Z}}c(n)|u^-(n)|^q\nonumber\\
&    & +\left(\frac{1-s^p}{p}-\frac{1-s^q}{q}\right)\sum_{n\in\mathbb{Z}}a(n)\Big|\sum_{i=1}^{\frac{p}{2}-1}\sum_{j=0}^i
             C_{\frac{p}{2}-1}^i C_i^j2^{i-j}(\Delta u^+(n))^{p-(i+j)}(\Delta u^-(n))^{i+j}\Big|\nonumber\\
&    & +\left(\frac{1-s^p}{p}-\frac{1-s^q}{q}\right)\sum_{n\in\mathbb{Z}}a(n)\Big|\sum_{i=1}^{\frac{p}{2}}\sum_{j=0}^{i-1}
             C_{\frac{p}{2}-1}^{i-1}C_{i-1}^{j}2^{i-1-j}(\Delta u^+(n))^{p-(i+j)}(\Delta u^-(n))^{i+j}\Big|\nonumber\\
&    & +\left(\frac{1-t^p}{p}-\frac{1-t^q}{q}\right)\sum_{n\in\mathbb{Z}}a(n)\Big|\sum_{i=1}^{\frac{p}{2}-1}\sum_{j=1}^{i}
             C_{\frac{p}{2}-1}^{i-1}C_{i-1}^{j-1}2^{i-j}(\Delta u^+(n))^{p-(i+j)}(\Delta u^-(n))^{i+j}\Big|\nonumber\\
&    & +\left(\frac{1-t^p}{p}-\frac{1-t^q}{q}\right)\sum_{n\in\mathbb{Z}}a(n)\Big|\sum_{j=1}^{\frac{p}{2}-1}
          C_{\frac{p}{2}-1}^{j-1}2^{\frac{p}{2}-j}(\Delta u^+(n))^{\frac{p}{2}-j}(\Delta u^-(n))^{\frac{p}{2}+j}\Big|\nonumber\\
&    & +\left(\frac{1-t^p}{p}-\frac{1-t^q}{q}\right)\sum_{n\in\mathbb{Z}}a(n)\Big|\sum_{i=1}^{\frac{p}{2}}\sum_{j=0}^{i-1}
             C_{\frac{p}{2}-1}^{i-1}C_{i-1}^{j}2^{i-1-j}(\Delta u^+(n))^{p-(i+j)}(\Delta u^-(n))^{i+j}\Big|\nonumber\\
&    &-\frac{1-s^p}{p}\sum_{n\in\mathbb{Z}}a(n)\Big|\sum_{i=1}^{\frac{p}{2}-1}\sum_{j=0}^i C_{\frac{p}{2}-1}^i C_i^j2^{i-j}
              (\Delta u^+(n))^{p-(i+j)}(\Delta u^-(n))^{i+j}\Big|\nonumber\\
&    & -\frac{1-s^p}{p}\sum_{n\in\mathbb{Z}}a(n)\Big|\sum_{i=1}^{\frac{p}{2}}\sum_{j=0}^{i-1} C_{\frac{p}{2}-1}^{i-1}
              C_{i-1}^{j}2^{i-1-j}(\Delta u^+(n))^{p-(i+j)}(\Delta u^-(n))^{i+j}\Big|\nonumber\\
&    & -\frac{1-t^p}{p}\sum_{n\in\mathbb{Z}}a(n)\Big|\sum_{i=1}^{\frac{p}{2}-1}\sum_{j=1}^{i}C_{\frac{p}{2}-1}^{i-1}
              C_{i-1}^{j-1}2^{i-j}(\Delta u^+(n))^{p-(i+j)}(\Delta u^-(n))^{i+j}\Big|\nonumber\\
&    & -\frac{1-t^p}{p}\sum_{n\in\mathbb{Z}}a(n)\Big|\sum_{j=1}^{\frac{p}{2}-1}C_{\frac{p}{2}-1}^{j-1}2^{\frac{p}{2}-j}
             (\Delta u^+(n))^{\frac{p}{2}-j}(\Delta u^-(n))^{\frac{p}{2}+j}\Big|\nonumber\\
&    & -\frac{1-t^p}{p}\sum_{n\in\mathbb{Z}}a(n)\Big|\sum_{i=1}^{\frac{p}{2}}\sum_{j=0}^{i-1} C_{\frac{p}{2}-1}^{i-1}
              C_{i-1}^{j}2^{i-1-j}(\Delta u^+(n))^{p-(i+j)}(\Delta u^-(n))^{i+j}\Big|\nonumber\\
&    & +\frac{1}{p}\sum_{n\in\mathbb{Z}}a(n)\Big|\sum_{i=1}^{\frac{p}{2}-1}\sum_{j=0}^iC_{\frac{p}{2}}^iC_i^j2^{i-j}
           (\Delta u^+(n))^{p-(i+j)}(\Delta u^-(n))^{i+j}\Big|\nonumber\\
&   &  +\frac{1}{p}\sum_{n\in\mathbb{Z}}a(n)\Big|\sum_{j=0}^{\frac{p}{2}-1}C_{\frac{p}{2}}^j2^{\frac{p}{2}-j}
           (\Delta u^+(n))^{\frac{p}{2}-j}(\Delta u^-(n))^{\frac{p}{2}+j}\Big|\nonumber\\
&   &  -\frac{1}{p}\sum_{n\in\mathbb{Z}}a(n)\Big|\sum_{i=1}^{\frac{p}{2}-1}\sum_{j=0}^iC_{\frac{p}{2}}^iC_i^j2^{i-j}
           (\Delta su^+(n))^{p-(i+j)}(\Delta tu^-(n))^{i+j}\Big|\nonumber\\
&   &  -\frac{1}{p}\sum_{n\in\mathbb{Z}}a(n)\Big|\sum_{j=0}^{\frac{p}{2}-1}C_{\frac{p}{2}}^j2^{\frac{p}{2}-j}
           (\Delta su^+(n))^{\frac{p}{2}-j}(\Delta tu^-(n))^{\frac{p}{2}+j}\Big|\nonumber\\
&\geq& \frac{1-s^{q}}{q}\langle I'(u),u^+\rangle+\frac{1-t^{q}}{q}\langle I'(u),u^-\rangle
       +\left(\frac{1-s^p}{p}-\frac{1-s^{q}}{q}\right)\|u^+\|^p+\left(\frac{1-t^p}{p}-\frac{1-t^{q}}{q}\right)\|u^-\|^p\nonumber\\
&    & -\frac{1-s^p}{p}\sum_{n\in\mathbb{Z}}a(n)\Big|\sum_{i=1}^{\frac{p}{2}-1}\sum_{j=0}^i C_{\frac{p}{2}-1}^i C_i^j2^{i-j}
              (\Delta u^+(n))^{p-(i+j)}(\Delta u^-(n))^{i+j}\Big|\nonumber\\
&    & -\frac{1-s^p}{p}\sum_{n\in\mathbb{Z}}a(n)\Big|\sum_{i=1}^{\frac{p}{2}}\sum_{j=0}^{i-1} C_{\frac{p}{2}-1}^{i-1}
              C_{i-1}^{j}2^{i-1-j}(\Delta u^+(n))^{p-(i+j)}(\Delta u^-(n))^{i+j}\Big|\nonumber\\
&    & -\frac{1-t^p}{p}\sum_{n\in\mathbb{Z}}a(n)\Big|\sum_{i=1}^{\frac{p}{2}-1}\sum_{j=1}^{i}C_{\frac{p}{2}-1}^{i-1}
              C_{i-1}^{j-1}2^{i-j}(\Delta u^+(n))^{p-(i+j)}(\Delta u^-(n))^{i+j}\Big|\nonumber\\
&    & -\frac{1-t^p}{p}\sum_{n\in\mathbb{Z}}a(n)\Big|\sum_{j=1}^{\frac{p}{2}-1}C_{\frac{p}{2}-1}^{j-1}2^{\frac{p}{2}-j}
             (\Delta u^+(n))^{\frac{p}{2}-j}(\Delta u^-(n))^{\frac{p}{2}+j}\Big|\nonumber\\
&    & -\frac{1-t^p}{p}\sum_{n\in\mathbb{Z}}a(n)\Big|\sum_{i=1}^{\frac{p}{2}}\sum_{j=0}^{i-1} C_{\frac{p}{2}-1}^{i-1}
              C_{i-1}^{j}2^{i-1-j}(\Delta u^+(n))^{p-(i+j)}(\Delta u^-(n))^{i+j}\Big|\nonumber\\
&    & +\frac{1}{p}\sum_{n\in\mathbb{Z}}a(n)\Big|\sum_{i=1}^{\frac{p}{2}-1}\sum_{j=0}^iC_{\frac{p}{2}}^iC_i^j2^{i-j}
           (\Delta u^+(n))^{p-(i+j)}(\Delta u^-(n))^{i+j}\Big|\nonumber\\
&   &  +\frac{1}{p}\sum_{n\in\mathbb{Z}}a(n)\Big|\sum_{j=0}^{\frac{p}{2}-1}C_{\frac{p}{2}}^j2^{\frac{p}{2}-j}
           (\Delta u^+(n))^{\frac{p}{2}-j}(\Delta u^-(n))^{\frac{p}{2}+j}\Big|\nonumber\\
&   &  -\frac{1}{p}\sum_{n\in\mathbb{Z}}a(n)\Big|\sum_{i=1}^{\frac{p}{2}-1}\sum_{j=0}^iC_{\frac{p}{2}}^iC_i^j2^{i-j}
           (\Delta su^+(n))^{p-(i+j)}(\Delta tu^-(n))^{i+j}\Big|\nonumber\\
&   &  -\frac{1}{p}\sum_{n\in\mathbb{Z}}a(n)\Big|\sum_{j=0}^{\frac{p}{2}-1}C_{\frac{p}{2}}^j2^{\frac{p}{2}-j}
           (\Delta su^+(n))^{\frac{p}{2}-j}(\Delta tu^-(n))^{\frac{p}{2}+j}\Big|\nonumber\\
& = & \frac{1-s^{q}}{q}\langle I'(u),u^+\rangle+\frac{1-t^{q}}{q}\langle I'(u),u^-\rangle
       +\left(\frac{1-s^p}{p}-\frac{1-s^{q}}{q}\right)\|u^+\|^p+\left(\frac{1-t^p}{p}-\frac{1-t^{q}}{q}\right)\|u^-\|^p\nonumber\\
&   & +\sum_{n\in\mathbb{Z}}a(n)\Big|\sum_{i=1}^{\frac{p}{2}-1}\sum_{j=1}^{i-1}2^{i-j}(\Delta u^+(n))^{p-(i+j)}
                   (\Delta u^-(n))^{i+j}\Big|\Theta'\nonumber\\
&   & +\frac{2C_{\frac{p}{2}}^{i}-2C_{\frac{p}{2}-1}^{i}-C_{\frac{p}{2}-1}^{i-1}-C_{\frac{p}{2}-1}^{i-1}}{2p}
            \sum_{n\in\mathbb{Z}}a(n)\Big|\sum_{i=1}^{\frac{p}{2}-1}2^i(\Delta u^+(n))^{p-i}(\Delta u^-(n))^i\Big|\nonumber\\
&   & +\frac{C_{\frac{p}{2}}^{i}-C_{\frac{p}{2}-1}^{i}-C_{\frac{p}{2}-1}^{i-1}}{p}
            \sum_{n\in\mathbb{Z}}a(n)\Big|\sum_{i=1}^{\frac{p}{2}-1}(\Delta u^+(n))^{p-2i}(\Delta u^-(n))^{2i}\Big|\nonumber\\
&   & +\frac{2C_{\frac{p}{2}}^{j}-C_{\frac{p}{2}-1}^{j}-2C_{\frac{p}{2}-1}^{j-1}-C_{\frac{p}{2}-1}^{j}}{2p}
             \sum_{n\in\mathbb{Z}}a(n)\Big|\sum_{j=0}^{\frac{p}{2}-1}2^{\frac{p}{2}-j}
            (\Delta u^+(n))^{\frac{p}{2}-j}(\Delta u^-(n))^{\frac{p}{2}+j}\Big|\nonumber\\
&   & +\sum_{n\in\mathbb{Z}}a(n)\Big|\sum_{i=1}^{\frac{p}{2}-1}\sum_{j=1}^{i-1}2^{i-j}(\Delta u^+(n))^{p-(i+j)}
                   (\Delta u^-(n))^{i+j}\Big|\Theta\nonumber\\
&   & +\frac{2s^pC_{\frac{p}{2}-1}^{i}+s^pC_{\frac{p}{2}-1}^{i-1}+t^pC_{\frac{p}{2}-1}^{i-1}-2s^{p-i}t^iC_{\frac{p}{2}}^{i}}{2p}
            \sum_{n\in\mathbb{Z}}a(n)\Big|\sum_{i=1}^{\frac{p}{2}-1}2^i(\Delta u^+(n))^{p-i}(\Delta u^-(n))^i\Big|\nonumber\\
&   & +\frac{s^p C_{\frac{p}{2}-1}^{i}+t^pC_{\frac{p}{2}-1}^{i-1}-s^{p-2i}t^{2i}C_{\frac{p}{2}}^{i}}{p}
            \sum_{n\in\mathbb{Z}}a(n)\Big|\sum_{i=1}^{\frac{p}{2}-1}(\Delta u^+(n))^{p-2i}(\Delta u^-(n))^{2i}\Big|\nonumber\\
&   & +\frac{s^pC_{\frac{p}{2}-1}^{j}+2t^pC_{\frac{p}{2}-1}^{j-1}+t^pC_{\frac{p}{2}-1}^{j}-2s^{\frac{p}{2}-j}t^{\frac{p}{2}+j}
             C_{\frac{p}{2}}^{j}}{2p}\sum_{n\in\mathbb{Z}}a(n)\Big|\sum_{j=0}^{\frac{p}{2}-1}2^{\frac{p}{2}-j}
            (\Delta u^+(n))^{\frac{p}{2}-j}(\Delta u^-(n))^{\frac{p}{2}+j}\Big|\nonumber\\
&\geq& \frac{1-s^{q}}{q}\langle I'(u),u^+\rangle+\frac{1-t^{q}}{q}\langle I'(u),u^-\rangle
       +\left(\frac{1-s^p}{p}-\frac{1-s^{q}}{q}\right)\|u^+\|^p+\left(\frac{1-t^p}{p}-\frac{1-t^{q}}{q}\right)\|u^-\|^p\nonumber\\
&   & +\sum_{n\in\mathbb{Z}}a(n)\Big|\sum_{i=1}^{\frac{p}{2}-1}\sum_{j=1}^{i-1}2^{i-j}(\Delta u^+(n))^{p-(i+j)}
                   (\Delta u^-(n))^{i+j}\Big|\Theta,
\end{eqnarray}
where $\Theta'=\frac{2C_{\frac{p}{2}}^i C_i^j-2C_{\frac{p}{2}-1}^i C_i^j-C_{\frac{p}{2}-1}^{i-1}C_{i-1}^j
-2C_{\frac{p}{2}-1}^{i-1}C_{i-1}^{j-1}-C_{\frac{p}{2}-1}^{i-1}C_{i-1}^j}{2p}=0$ (use combination formula) and $\Theta\ge0$ (see Appendix 2). So, we obtain that (\ref{Eq5}) holds for all $u\in\mathcal{D},s,t\geq0$.
\qed
\vskip2mm
\noindent
{\bf Remark 2.1.} {\it Let $s=t$ in} (\ref{Eq5}). {\it It is easy to see that $\Theta=0$. Then for all $u\in\mathcal{D}$ and $t\geq0$, there holds}
\begin{eqnarray*}
        I(u)
& \geq & I(tu)+\frac{1-t^{q}}{q}\langle I'(u),u\rangle+\left(\frac{1-t^p}{p}-\frac{1-t^{q}}{q}\right)\|u^+\|^p
           +\left(\frac{1-t^p}{p}-\frac{1-t^{q}}{q}\right)\|u^-\|^p.
\end{eqnarray*}

\vskip2mm
\noindent
{\bf Corollary 2.3.} {\it Assume that $(C_1)$ and $(C_2)$ hold. For all  $u\in\mathcal{D}$ and $t\geq0$, we have}
\begin{eqnarray}
\label{Eq9}
        I(u)
\geq  I(tu)+\frac{1-t^{q}}{q}\langle I'(u),u\rangle+\left(\frac{1-t^p}{p}-\frac{1-t^{q}}{q}\right)\|u\|^p.
\end{eqnarray}
{\bf Proof.}
According to (\ref{eq7}),  (\ref{eq10}) and (\ref{Eq6}),  there exists
\begin{eqnarray}
\label{Eq11}
&    & I(u)-I(tu)\nonumber\\
& =  & \frac{1}{p}\left(\|u\|^p-\|tu\|^p\right)
         +\frac{r}{q^2}\sum_{n\in\mathbb{Z}}c(n)\left(|u(n)|^q-|tu(n)|^q\right)\nonumber\\
&    &  -\frac{1}{q}\sum_{n\in\mathbb{Z}}c(n)\left(|u(n)|^q\ln{|u(n)|^r-|tu(n)|^q\ln{|tu(n)|^r}}\right)\nonumber\\
& =  & \frac{1-t^p}{p}\|u\|^p+\frac{r(1-t^q)}{q^2}\sum_{n\in\mathbb{Z}}c(n)|u(n)|^q
          -\frac{1-t^q}{q}\sum_{n\in\mathbb{Z}}c(n)|u|^q\ln{|u|^r}+\frac{1}{q}\sum_{n\in\mathbb{Z}}c(n)|tu|^q\ln{t^r}\nonumber\\
& =  & \frac{1-t^q}{q}\big(\langle I'(u),u\rangle-\|u\|^p\big)+\frac{1-t^p}{p}\|u\|^p
         +\frac{r(1-t^q)}{q^2}\sum_{n\in\mathbb{Z}}c(n)|u(n)|^q+\frac{1}{q}\sum_{n\in\mathbb{Z}}c(n)|tu|^q\ln{t^r}\nonumber\\
& =  & \frac{1-t^q}{q}\langle I'(u),u\rangle+\left(\frac{1-t^p}{p}-\frac{1-t^q}{q}\right)\|u\|^p
           +\frac{r(1-t^q)+qt^q\ln t^r}{q^2}\sum_{n\in\mathbb{Z}}c(n)|u(n)|^q\nonumber\\
&\geq& \frac{1-t^q}{q}\langle I'(u),u\rangle+\left(\frac{1-t^p}{p}-\frac{1-t^q}{q}\right)\|u\|^p
\end{eqnarray}
Hence, (\ref{Eq9}) holds for all $u\in\mathcal{D}$ and  $t\geq0$.
\qed
\par
Note that $1<p<q$, $\Theta\geq0$ and the function $f(x)=\frac{1-a^x}{x}$ is strictly monotonically decreasing on $(0,+\infty)$ for $a>0$ and $a\neq 1$. Then in combination with Lemma 2.2, we have the following corollary.
\vskip2mm
\noindent
{\bf Corollary 2.4.} {\it Assume that $(C_1)$ and $(C_2)$ hold. For any $u\in\mathcal{M}$, we can obtain that $I(u)=\max_{s,t\geq0}I(su^++tu^-)$.}
\par
In combination with Corollary 2.3 or Remark 2.1, we have the following corollary.
\vskip2mm
\noindent
{\bf Corollary 2.5.} {\it Assume that $(C_1)$ and $(C_2)$ hold. For any $u\in\mathcal{N}$, there holds $I(u)=\max_{t\geq0}I(tu)$.}
\vskip2mm
\noindent
{\bf Lemma 2.6.} {\it Assume that $(C_1)$ and $(C_2)$ hold. For any $u\in\mathcal{D}$ with $u\neq0$, there exists a unique positive constant $t_0$ such that $t_0u\in\mathcal{N}$.}\\
{\bf Proof.} First, we prove that the existence of $t_0$. For any $u\in\mathcal{D}$ with $u\neq0$, let $u\in\mathcal{N}$ be fixed and define a function $g(t)=\langle I'(tu),tu \rangle$ on $(0,+\infty)$. On one hand, by (\ref{eq9}) and Lemma 2.1, there exist two positive contants $\varepsilon_2<\frac{b_0}{c_0}$ and $C_{\varepsilon_2}$ such that
\begin{eqnarray}
\label{Eq12}
       g(t)
& =  & t^p\|u\|^p-\sum_{n\in\mathbb{Z}}c(n)|tu(n)|^q\ln{|tu(n)|^r}\nonumber\\
&\geq& t^p\|u\|^p-\sum_{n\in\mathbb{Z}}c_0\varepsilon_2|tu(n)|^p
                             -\sum_{n\in\mathbb{Z}}c_0C_{\varepsilon_2}|tu(n)|^\zeta\nonumber\\
&\geq& t^p\|u\|^p-b_0^{-1}c_0\varepsilon_2 t^p\|u\|^p-t^\zeta\sum_{n\in\mathbb{Z}}c_0C_{\varepsilon_2}|u(n)|^\zeta.
\end{eqnarray}
Then, according to $\zeta>q$ and $q>p>1$, we have that $g(t)>0$ for all sufficiently small $t>0$.
\par
On the other hand, noting that $c(n)>0$ for all $n\in\mathbb{Z}$, by $(C_2)$ and (\ref{eq9}), there exists
\begin{eqnarray}
\label{Eq13}
       g(t)
& =  & t^p\|u\|^p-t^q\ln t^r\sum_{n\in\mathbb{Z}}c(n)|u(n)|^q-t^q\sum_{n\in\mathbb{Z}}c(n)|u(n)|^q\ln{|u(n)|^r}\nonumber\\
&\leq& t^p\|u\|^p-t^q\ln t^r\sum_{n\in\mathbb{Z}}c(n)|u(n)|^q
           +t^q\sum_{\substack{n\in\mathbb{Z}}}c_0(\varepsilon|u(n)|^p+C_\varepsilon|u(n)|^\zeta).
\end{eqnarray}
Then, by $1<p<q$, $r\ge 1$ and (\ref{Eq13}), it is easy to see $g(t)<0$ for all large $t$. So it follows from the continuity of $g(t)$ that there exists a $t_0\in(0,+\infty)$ such that $g(t_0)=0$, which implies that there exists a positive constant $t_0$ such that $t_0u\in\mathcal{N}$.
\par
Next, we prove the uniqueness of $t_0$. Proofing by contradiction, we assume that there exist $u\in\mathcal{D}$ and two positive numbers $t_1\neq t_2$ such that $t_1u\in\mathcal{N}$ and $t_2u\in\mathcal{N}$. Note that the function $f(x)=\frac{1-a^x}{x}$ is strictly monotonically decreasing on $(0,+\infty)$ for $a>0$ and $a\neq1$. Taking $u$ as $t_1u$ and $t$ as $\frac{t_2}{t_1}$ in Corollary 2.3, there holds
\begin{eqnarray}
\label{Eq14}
      I(t_1u)
 \geq I(t_2u)+t_1^p\left(\frac{1-(\frac{t_2}{t_1})^p}{p}-\frac{1-(\frac{t_2}{t_1})^{q}}{q}\right)\|u\|^p
   >  I(t_2u).
\end{eqnarray}
On the other hand, taking $u$ as $t_2u$ and $t$ as $\frac{t_1}{t_2}$ in Corollary 2.3, there also holds
\begin{eqnarray}
\label{Eq15}
      I(t_2u)
 \geq I(t_1u)+t_2^p\left(\frac{1-(\frac{t_1}{t_2})^p}{p}-\frac{1-(\frac{t_1}{t_2})^{q}}{q}\right)\|u\|^p
   >  I(t_1u).
\end{eqnarray}
Hence, (\ref{Eq14}) contradicts with (\ref{Eq15}). So $t_1=t_2$, that is, there exists a unique positive constant $t_0$ such that $t_0u\in\mathcal{N}$.
\qed

\vskip2mm
\noindent
{\bf Lemma 2.7.} {\it Assume that $(C_1)$ and $(C_2)$ hold. For any $u\in\mathcal{D}$ with $u^\pm\neq 0$, there exists a unique pair of positive constants $(s_0,t_0)$ such that $s_0u^++t_0u^-\in\mathcal{M}$.}\\
{\bf Proof.} First, we prove that the existence of $(s_0,t_0)$. For any $u\in\mathcal{D}$ with $u^\pm\neq 0$, according to (\ref{Eq3}) and (\ref{Eq4}), we have
\begin{eqnarray}
\label{Eq16}
     h_1(s,t)
& := &\langle I'(su^++tu^-),su^+\rangle\nonumber\\
& = &s^p\|u^+\|^p-\sum_{n\in\mathbb{Z}}c(n)|su^+(n)|^q\ln{|su^+(n)|^r}\nonumber\\
&   &+\sum_{n\in\mathbb{Z}}a(n)\Big|\sum_{i=1}^{\frac{p}{2}-1}\sum_{j=0}^i C_{\frac{p}{2}-1}^i C_i^j2^{i-j}
                 s^{p-(i+j)}t^{i+j}(\Delta u^+(n))^{p-(i+j)}(\Delta u^-(n))^{i+j}\Big|\nonumber\\
&  &   +\sum_{n\in\mathbb{Z}}a(n)\Big|\sum_{i=1}^{\frac{p}{2}}\sum_{j=0}^{i-1} C_{\frac{p}{2}-1}^{i-1} C_{i-1}^{j}2^{i-1-j}
                 s^{p-(i+j)}t^{i+j}(\Delta u^+(n))^{p-(i+j)}(\Delta u^-(n))^{i+j}\Big|,
\end{eqnarray}
and
\begin{eqnarray}
\label{Eq17}
     h_2(s,t)
& := &\langle I'(su^++tu^-),tu^-\rangle\nonumber\\
& = &t^p\|u^-\|^p-\sum_{n\in\mathbb{Z}}c(n)|tu^-(n)|^q\ln{|tu^-(n)|^r}\nonumber\\
&   &+\sum_{n\in\mathbb{Z}}a(n)\Big|\sum_{i=1}^{\frac{p}{2}-1}\sum_{j=1}^{i} C_{\frac{p}{2}-1}^{i-1}
               C_{i-1}^{j-1}2^{i-j}s^{p-(i+j)}t^{i+j}(\Delta u^+(n))^{p-(i+j)}(\Delta u^-(n))^{i+j}\Big|\nonumber\\
&  &   +\sum_{n\in\mathbb{Z}}a(n)\Big|\sum_{j=1}^{\frac{p}{2}-1}C_{\frac{p}{2}-1}^{j-1}2^{\frac{p}{2}-j}
               s^{\frac{p}{2}-j}t^{\frac{p}{2}+j}(\Delta u^+(n))^{\frac{p}{2}-j}(\Delta u^-(n))^{\frac{p}{2}+j}\Big|\nonumber\\
&  &   +\sum_{n\in\mathbb{Z}}a(n)\Big|\sum_{i=1}^{\frac{p}{2}}\sum_{j=0}^{i-1} C_{\frac{p}{2}-1}^{i-1} C_{i-1}^{j}2^{i-1-j}
               s^{p-(i+j)}t^{i+j}(\Delta u^+(n))^{p-(i+j)}(\Delta u^-(n))^{i+j}\Big|.
\end{eqnarray}
It follows from (\ref{Eq12}) and (\ref{Eq13}) that $h_1(s,s)>0$ and $h_2(s,s)>0$ for $s>0$ sufficiently small and $h_1(t,t)<0$ and $h_2(t,t)<0$ for $t>0$ large enough. Thus, there are two constants $0<\theta_1<\theta_2$ such that
\begin{eqnarray}
\label{Eq18}
h_1(\theta_1,\theta_1)>0,\;h_2(\theta_1,\theta_1)>0, \;h_1(\theta_2,\theta_2)<0,\;h_2(\theta_2,\theta_2)<0.
\end{eqnarray}
For all $s,t\in[\theta_1,\theta_2]$, according to (\ref{Eq16}), (\ref{Eq17}) and (\ref{Eq18}), there exists
\begin{eqnarray}
\label{Q1}
      h_1(\theta_1,t)
&\geq& \theta_1^p\|u^+\|^p-\sum_{n\in\mathbb{Z}}c(n)|\theta_1u^+(n)|^q\ln{|\theta_1u^+(n)|^r}\nonumber\\
&   & +\sum_{n\in\mathbb{Z}}a(n)\Big|\sum_{i=1}^{\frac{p}{2}-1}\sum_{j=0}^i C_{\frac{p}{2}-1}^i C_i^j2^{i-j}
                 \theta_1^{p}(\Delta u^+(n))^{p-(i+j)}(\Delta u^-(n))^{i+j}\Big|\nonumber\\
&  &   +\sum_{n\in\mathbb{Z}}a(n)\Big|\sum_{i=1}^{\frac{p}{2}}\sum_{j=0}^{i-1} C_{\frac{p}{2}-1}^{i-1} C_{i-1}^{j}2^{i-1-j}
                 \theta_1^{p}(\Delta u^+(n))^{p-(i+j)}(\Delta u^-(n))^{i+j}\Big|\nonumber\\
& = & h_1(\theta_1,\theta_1)\nonumber\\
& > & 0,
\end{eqnarray}
and similarly, we can obtain that
\begin{eqnarray}
\label{Q2}
h_1(\theta_2,t)\leq h_1(\theta_2,\theta_2)<0,\;\;h_2(s,\theta_1)\geq h_2(\theta_1,\theta_1)>0\;\;\text{and}\;\;h_2(s,\theta_2)\leq h_2(\theta_2,\theta_2)<0.
\end{eqnarray}
Therefore, in virtue of (\ref{Q1}), (\ref{Q2}) and the Pincare-Miranda Theorem \cite{Miranda 1940}, there appears a point $(s_0,t_0)$ with $\theta_1<s_0,t_0<\theta_2$ such that $h_1(s_0,t_0)=h_2(s_0,t_0)=0$, that is, there exists a pair of positive constants $(s_0,t_0)$ such that $s_0u^++t_0u^-\in\mathcal{M}$.
\par
Next, we prove the uniquness of $(s_0,t_0)$. Proofing by contradiction, we suppose that there are two unequal pairs of positive constants $(s_1,t_1)$ and $(s_2,t_2)$ such that $s_1u^++t_1u^-\in\mathcal{M}$ and $s_2u^++t_2u^-\in\mathcal{M}$. Note that the function $f(x)=\frac{1-a^x}{x}$ is strictly monotonically decreasing on $(0,+\infty)$ for $a>0$ and $a\neq1$. For one thing, taking $u$, $s$ and $t$ as $s_1u^++t_1u^-$, $\frac{s_2}{s_1}$ and $\frac{t_2}{t_1}$ in Lemma 2.2, respectively, and noting that $p<q$, then we have
\begin{small}
\begin{eqnarray}
\label{Q5}
         I(s_1u^++t_1u^-)
&\geq& I(s_2u^++t_2u^-)+s_1^p\left(\frac{1-(\frac{s_2}{s_1})^p}{p}-\frac{1-(\frac{s_2}{s_1})^{q}}{q}\right)\|u^+\|^p
          +t_1^p\left(\frac{1-(\frac{t_2}{t_1})^p}{p}-\frac{1-(\frac{t_2}{t_1})^{q}}{q}\right)\|u^-\|^p\nonumber\\
&    & +\sum_{n\in\mathbb{Z}}a(n)\Big|\sum_{i=1}^{\frac{p}{2}-1}\sum_{j=1}^{i-1}2^{i-j}(\Delta u^+(n))^{p-(i+j)}
              (\Delta u^-(n))^{i+j}\Big|s_1^{p-{i+j}}t_1^{i+j}\Theta''\nonumber\\
&  > & I(s_2u^++t_2u^-),
\end{eqnarray}
\end{small}where $\Theta''=\frac{2(\frac{s_2}{s_1})^pC_{\frac{p}{2}-1}^iC_i^j+(\frac{s_2}{s_1})^pC_{\frac{p}{2}-1}^{i-1}C_{i-1}^j
+2(\frac{t_2}{t_1})^pC_{\frac{p}{2}-1}^{i-1}C_{i-1}^{j-1}+(\frac{t_2}{t_1})^pC_{\frac{p}{2}-1}^{i-1}C_{i-1}^j
-2(\frac{s_2}{s_1})^{p-(i+j)}(\frac{t_2}{t_1})^{i+j}C_{\frac{p}{2}}^iC_i^j}{2p}\geq0$ (see Appendix 3).
For another, taking $u$, $s$ and $t$ as $s_2u^++t_2u^-$, $\frac{s_1}{s_2}$ and $\frac{t_1}{t_2}$, respectively,  we also have
\begin{eqnarray}
\label{Q6}
I(s_2u^++t_2u^-)>I(s_1u^++t_1u^-).
\end{eqnarray}
As a consequence, there is a contradiction between $(\ref{Q5})$ and $(\ref{Q6})$. So $(s_1,t_1)=(s_2,t_2)$ which implies that there is a unique pair of positive constants $(s_0,t_0)$ such that $s_0u^++t_0u^-\in\mathcal{M}$.
\qed
\vskip2mm
\noindent
{\bf Lemma 2.8.} {\it Assume that $(C_1)$ and $(C_2)$ hold. Then}
\begin{eqnarray}
\label{W1}
  \inf_{u\in\mathcal{N}}I(u)
= :c_*
= \inf_{u\in \mathcal{D},u\neq 0}\max_{t\geq 0}I(tu).
\end{eqnarray}
and
\begin{eqnarray}
\label{Q7}
  \inf_{u\in\mathcal{M}}I(u)
= :m_*
= \inf_{u\in \mathcal{D},u^\pm\neq 0}\max_{s,t\geq 0}I(su^++tu^-).
\end{eqnarray}
\vskip0mm
\noindent
{\bf Proof.} On one hand, according to Corollary 2.4 and the definition of $\mathcal{M}$, there holds
$$
     \inf_{u\in\mathcal{M}}I(u)
  =  \inf_{u\in \mathcal{M}} \max_{s,t\geq 0}I(su^++tu^-)
\geq \inf_{u\in \mathcal{D},u^\pm\neq 0}\max_{s,t\geq 0}I(su^++tu^-).
$$
On the other hand, for any $u\in \mathcal{D}$ with $u^\pm\neq0$, it follows from Lemma 2.7 that there appears two positive constants $s_0,t_0$ such that $s_0u^++t_0u^-\in\mathcal{M}$. Then we have
$$
     \max_{s,t\ge 0} I(su^++tu^-)
\geq I(s_0 u^++t_0 u^-)
\geq \inf_{u\in\mathcal{M}}I(u),
$$
which implies that
$$
     \inf_{u\in \mathcal{D},u^\pm\neq0}\max_{s,t\ge 0} I(su^++tu^-)
\geq \inf_{u\in \mathcal{D},u^\pm\neq0}I(s_0 u^++t_0 u^-)
\geq \inf_{u\in\mathcal{M}}I(u).
$$
Hence, it is easy to see that the conclusion $(\ref{Q7})$ holds. Similarly, it follows from Corollary 2.5, the definition of $\mathcal{N}$ and Lemma 2.6 that $(\ref{W1})$ also holds.
\qed

\vskip2mm
\noindent
{\bf Lemma 2.9.} {\it Assume that $(C_1)$ and $(C_2)$ hold. Then $m_*>0$ and $c_*>0$ can be achieved.}
\vskip0mm
\noindent
{\bf Proof.} For any $u\in\mathcal{M}$, there holds $\langle I'(u),u\rangle=0$. For $\varepsilon_3=\frac{b_0}{pc_0}>0$, it follows from (\ref{eq9}), (\ref{eq10}) and  Lemma 2.1 that there is a positive constant $C_{\varepsilon_3}$ such that
\begin{eqnarray*}
\label{Q8}
\|u\|^p=\sum_{n\in\mathbb{Z}}c(n)|u(n)|^{q}\ln |u(n)|^r
       \leq c_0\varepsilon_3\|u\|_{l^p}^{p}+c_0C_{\varepsilon_3}\|u\|_{l^\zeta}^{\zeta}
       \leq \frac{1}{p}\|u\|^p+c_0C_{\varepsilon_3} b_0^{-\frac{\zeta}{p}}\|u\|^\zeta.
\end{eqnarray*}
Since $1<p<q<\zeta$, then $\|u\|\geq \rho:=\left(\frac{(p-1)b_0^{\frac{\zeta}{p}}}{pc_0C_{\varepsilon_3}}\right)^{\zeta-p}$ for any $u\in\mathcal{M}$.
\par
Let $\{u_k\}\subset \mathcal{M}$ be such that $I(u_k)\rightarrow m_*.$ By (\ref{eq7}) and (\ref{eq10}), there holds
\begin{eqnarray*}
m_*+o(1)=I(u_k)-\frac{1}{q}\langle I'(u_k),u_k\rangle
        =\left(\frac{1}{p}-\frac{1}{q}\right)\|u_k\|^p+\frac{r}{q^2}\sum_{n\in\mathbb{Z}}c(n)|u_k(n)|^q
        \geq\left(\frac{1}{p}-\frac{1}{q}\right)\|u_k\|^p.
\end{eqnarray*}
This shows that the sequence $\{u_k\}$ is bounded in $\mathcal{D}$, that is, there exists a $M_1>0$ such that $\|u_k\|\leq M_1$. Thus, there appears a $u_0\in\mathcal{D}$ such that $u_k^\pm\rightharpoonup u_0^\pm$ in $\mathcal{D}$. Then it follows from Lemma 2.1 that $u_k^\pm\rightarrow u_0^\pm$ in $l^\kappa(\mathbb{Z},\mathbb{R})$ for $\kappa\in[p,+\infty]$ and $u_k^\pm(n)\rightarrow u_0^\pm(n)$ for all $n\in\mathbb{Z}$.
\par
Since $\{u_k\}\subset \mathcal{M}$, there exists $\langle I'(u_k),u_k^\pm\rangle=0$ and then by Proposition 2.1, we have
\begin{eqnarray}
\label{Q9}
&   &      \|u_k^+\|^p-\sum_{n\in\mathbb{Z}}c(n)|u^+_k(n)|^q\ln{|u^+_k(n)|^r}\nonumber\\
& = & -\sum_{n\in\mathbb{Z}}a(n)\Big|\sum_{i=1}^{\frac{p}{2}-1}\sum_{j=0}^i C_{\frac{p}{2}-1}^i C_i^j2^{i-j}
              (\Delta u_k^+(n))^{p-(i+j)}(\Delta u_k^-(n))^{i+j}\Big|\nonumber\\
&   & -\sum_{n\in\mathbb{Z}}a(n)\Big|\sum_{i=1}^{\frac{p}{2}}\sum_{j=0}^{i-1} C_{\frac{p}{2}-1}^{i-1} C_{i-1}^{j}2^{i-1-j}
              (\Delta u_k^+(n))^{p-(i+j)}(\Delta u_k^-(n))^{i+j}\Big|
\end{eqnarray}
and
\begin{eqnarray}
\label{Q10}
 & &       \|u_k^-\|^p-\sum_{n\in\mathbb{Z}}c(n)|u^-_k(n)|^q\ln{|u^-_k(n)|^r}\nonumber\\
& = &   -\sum_{n\in\mathbb{Z}}a(n)\Big|\sum_{i=1}^{\frac{p}{2}-1}\sum_{j=1}^{i} C_{\frac{p}{2}-1}^{i-1}
                       C_{i-1}^{j-1}2^{i-j}(\Delta u_k^+(n))^{p-(i+j)}(\Delta u_k^-(n))^{i+j}\Big|\nonumber\\
&   &   -\sum_{n\in\mathbb{Z}}a(n)\Big|\sum_{j=1}^{\frac{p}{2}-1}C_{\frac{p}{2}-1}^{j-1}2^{\frac{p}{2}-j}
                       (\Delta u_k^+(n))^{\frac{p}{2}-j}(\Delta u_k^-(n))^{\frac{p}{2}+j}\Big|\nonumber\\
&   &   -\sum_{n\in\mathbb{Z}}a(n)\Big|\sum_{i=1}^{\frac{p}{2}}\sum_{j=0}^{i-1} C_{\frac{p}{2}-1}^{i-1} C_{i-1}^{j}2^{i-1-j}
              (\Delta u_k^+(n))^{p-(i+j)}(\Delta u_k^-(n))^{i+j}\Big|.
\end{eqnarray}
It follows from (\ref{eq9}), (\ref{Q9}), Lemma 2.1 and the boundedness of $\{u_k\}$ that there exists $\varepsilon_4\in (0,\frac{b_0\rho^p}{c_0M_1^p})$ and a positive constant $C_{\varepsilon_4}$ such that
\begin{eqnarray*}
         \rho^p
& \leq &  \|u_k^+\|^p \nonumber\\
&  =   &  \sum_{n\in\mathbb{Z}}c(n)|u^+_k(n)|^q\ln{|u^+_k(n)|^r}
          -\sum_{n\in\mathbb{Z}}a(n)\Big|\sum_{i=1}^{\frac{p}{2}-1}\sum_{j=0}^i C_{\frac{p}{2}-1}^i C_i^j2^{i-j}
              (\Delta u_k^+(n))^{p-(i+j)}(\Delta u_k^-(n))^{i+j}\Big|\nonumber\\
&      &  -\sum_{n\in\mathbb{Z}}a(n)\Big|\sum_{i=1}^{\frac{p}{2}}\sum_{j=0}^{i-1} C_{\frac{p}{2}-1}^{i-1} C_{i-1}^{j}2^{i-1-j}
              (\Delta u_k^+(n))^{p-(i+j)}(\Delta u_k^-(n))^{i+j}\Big|\nonumber\\
& \leq &  \sum_{n\in\mathbb{Z}}c(n)|u^+_k(n)|^q\ln{|u^+_k(n)|^r}\nonumber\\
& \leq &  c_0\varepsilon_4\|u_k^+\|_{l^p}^{p}+c_0C_{\varepsilon_4}\|u_k^+\|_{l^\zeta}^{\zeta}\nonumber\\
& \leq &  c_0\varepsilon_4 b_0^{-1}\|u_k^+\|^{p}+c_0C_{\varepsilon_4}\|u_k^+\|_{l^\zeta}^{\zeta}\nonumber\\
& \leq &  c_0\varepsilon_4 b_0^{-1}M_1^{p}+c_0C_{\varepsilon_4}\|u_k^+\|_{l^\zeta}^{\zeta},
\end{eqnarray*}
which implies that $\|u_k^+\|_{l^\zeta}^{\zeta}\geq \frac{\rho^p-c_0\varepsilon_4b_0^{-1}M_1^p}{c_0C_{\varepsilon_4}}>0$. Similarly, by (\ref{eq9}), (\ref{Q10}), Lemma 2.1 and the boundedness of $\{u_k\}$, there exists $\varepsilon_5\in (0,\frac{b_0\rho^p}{c_0M_1^p})$ and a positive constant $C_{\varepsilon_5}$ such that $\|u_k^-\|_{l^\zeta}^{\zeta}\geq \frac{\rho^p-c_0\varepsilon_5 b_0^{-1}M_1^p}{c_0C_{\varepsilon_5}}>0$. Then let $\varepsilon':=\max\{\varepsilon_4,\varepsilon_5\}$ and $C_{\varepsilon'}:=\max\{C_{\varepsilon_4},C_{\varepsilon_5}\}$, we have that $\|u_k^\pm\|_{l^\zeta}^{\zeta}\geq \frac{\rho^p-c_0\varepsilon' b_0^{-1}M_1^p}{c_0C_{\varepsilon'}}>0$. For any $p\leq\kappa\leq+\infty$, it follows from the compactness of the embedding $\mathcal{D}\hookrightarrow l^\kappa(\mathbb{Z},\mathbb{R})$ that
$$
\|u_0^\pm\|_{l^\zeta}^{\zeta}\geq \frac{\rho^p-c_0\varepsilon' b_0^{-1}M_1^p}{c_0C_{\varepsilon'}}>0,
$$
which implies that $u_0^\pm\neq 0$. Note that $\Delta u^\pm(n)=u^\pm(n+1)-u^\pm(n)$. By the fact that $u_k^\pm(n)\rightarrow u_0^\pm(n)$ for all $n\in\mathbb{Z}$, we can derive that
\begin{eqnarray*}
    |\Delta u^\pm_k(n)-\Delta u^\pm_0(n)|= |u^\pm_k(n+1)-u^\pm_0(n+1)-(u^\pm_k(n)-u^\pm_0(n))|\to0,
\end{eqnarray*}
which implies that $\Delta u^\pm_k(n)\rightarrow\Delta u^\pm_0(n)$ for all $n\in\mathbb{Z}$. Note that
 \begin{eqnarray}\label{p1}
|u_k^+(n)|\leq\|u_k\|_{l^\infty}\leq b_0^{-\frac{1}{p}}\|u\|\leq b_0^{-\frac{1}{p}} M_1, \mbox{ for all } n\in\mathbb{Z}.
\end{eqnarray}
And by (\ref{eq9}), for any given $\varepsilon >0,$ there exists $C_\varepsilon>0$ such that
\begin{eqnarray}
\label{Q111}
|t|^q|\ln t^r|\leq\varepsilon|t|^p+C_\varepsilon|t|^\zeta, \forall t\in\mathbb{R}.
\end{eqnarray}
Then we can obtain that $|u_k^+(n)|^q|\ln |u_k^+(n)|^r|\leq\varepsilon |u_k^+(n)|^p+C_\varepsilon |u_k^+(n)|^\zeta\leq\varepsilon b_0^{-1}M_1^p+C_\varepsilon b_0^{-\frac{\zeta}{p}}M_1^\zeta$, which implies that
\begin{eqnarray*}
     \sum_{n\in\mathbb{Z}}c(n)|u_k^+(n)|^q|\ln |u_k^+(n)|^r|
\leq \left(\varepsilon b_0^{-1}M_1^p+C_\varepsilon b_0^{-\frac{\zeta}{p}}M_1^\zeta\right)\sum_{n\in\mathbb{Z}}c(n).
\end{eqnarray*}
Note that $\sum_{n\in\mathbb{Z}}c(n)<\infty$ (by $(C_2)$). Thus, it follows from $(\ref{Q9})$, the weak lower semi-continuity of norm, Fatou's Lemma and Lebesgue dominated convergence theorem that
\begin{eqnarray*}
&    & \langle I'(u_0),u_0^+\rangle\\
&\leq&\liminf_{k\rightarrow\infty}\|u_k^+\|^p
         +\liminf_{k\rightarrow\infty}\sum_{n\in\mathbb{Z}}a(n)\Big|\sum_{i=1}^{\frac{p}{2}-1}\sum_{j=0}^i
              C_{\frac{p}{2}-1}^i C_i^j2^{i-j}(\Delta u_k^+(n))^{p-(i+j)}(\Delta u_k^-(n))^{i+j}\Big|\\
&    &   +\liminf_{k\rightarrow\infty}\sum_{n\in\mathbb{Z}}a(n)\Big|\sum_{i=1}^{\frac{p}{2}}\sum_{j=0}^{i-1}
              C_{\frac{p}{2}-1}^{i-1} C_{i-1}^{j}2^{i-1-j}(\Delta u_k^+(n))^{p-(i+j)}(\Delta u_k^-(n))^{i+j}\Big|\\
&    &   -\liminf_{k\rightarrow\infty}\sum_{n\in\mathbb{Z}}c(n)|u^+_k(n)|^q\ln{|u^+_k(n)|^r}\\
&\leq&\liminf_{k\rightarrow\infty}\left[\|u_k^+\|^p
            +\sum_{n\in\mathbb{Z}}a(n)\Big|\sum_{i=1}^{\frac{p}{2}-1}\sum_{j=0}^i C_{\frac{p}{2}-1}^i C_i^j2^{i-j}
                   (\Delta u_k^+(n))^{p-(i+j)}(\Delta u_k^-(n))^{i+j}\Big| \right. \\
&    &\left.+\sum_{n\in\mathbb{Z}}a(n)\Big|\sum_{i=1}^{\frac{p}{2}}\sum_{j=0}^{i-1} C_{\frac{p}{2}-1}^{i-1} C_{i-1}^{j}2^{i-1-j}
              (\Delta u_k^+(n))^{p-(i+j)}(\Delta u_k^-(n))^{i+j}\Big|
              -\sum_{n\in\mathbb{Z}}c(n)|u^+_k(n)|^q\ln{|u^+_k(n)|^r}\right]\\
&  = & \liminf_{k\rightarrow\infty}\langle I'(u_k),u_k^+\rangle=0,
\end{eqnarray*}
which implies that
\begin{eqnarray}
\label{Q11}
\langle I'(u_0),u_0^+\rangle\leq 0.
\end{eqnarray}
Similarly, by $(\ref{Q10})$, the weak lower semi-continuity of norm, Fatou's Lemma and Lebesgue dominated convergence theorem, there exists
\begin{eqnarray}
\label{Q12}
\langle I'(u_0),u_0^-\rangle\leq 0\;\;\text{and then}\;\; \langle I'(u_0),u_0\rangle\leq 0.
\end{eqnarray}
According to Lemma 2.7, there are two positive constants $s_3,t_3$ such that
\begin{eqnarray}
\label{Q13}
s_3u_0^++t_3u_0^-\in\mathcal{M}\;\;\text{and}\;\;I(s_3u_0^++t_3u_0^-)\geq m_*.
\end{eqnarray}
By $(\ref{eq7})$, $(\ref{eq10})$, the weakly lower semi-continuity of norm, $(C_2)$, (\ref{p1}), Lemma 2.2, $(\ref{Q11})$, $(\ref{Q12})$ and $(\ref{Q13})$, there exists
\begin{eqnarray}
\label{Q14}
         m_*
&   =   &\lim_{k\rightarrow\infty}\left[I(u_k)-\frac{1}{q}\langle I'(u_k),u_k\rangle\right]\nonumber\\
& \geq  &\liminf_{k\rightarrow\infty}\left[\left(\frac{1}{p}-\frac{1}{q}\right)\|u_k\|^p
         +\frac{r}{q^2}\sum_{n\in\mathbb{Z}}c(n)|u_k(n)|^q\right]\nonumber\\
& \geq  &\left(\frac{1}{p}-\frac{1}{q}\right)\|u_0\|^p
         +\frac{r}{q^2}\sum_{n\in\mathbb{Z}}c(n)|u_0(n)|^q\nonumber\\
&   =   &I(u_0)-\frac{1}{q}\langle I'(u_0),u_0\rangle\nonumber\\
&  \geq &I(s_3u_0^++t_3u_0^-)+\frac{1-s_3^{q}}{q}\langle I'(u_0),u_0^+\rangle
        +\frac{1-t_3^{q}}{q}\langle I'(u_0),u_0^-\rangle
        -\frac{1}{q}\langle I'(u_0),u_0\rangle\nonumber\\
&  \geq &m_*-\frac{s_3^q}{q}\langle I'(u_0),u_0^+\rangle
        -\frac{t_3^q}{q}\langle I'(u_0),u_0^-\rangle\nonumber\\
&  \geq &m_*.
\end{eqnarray}
Moreover, in combination (\ref{Q14}) with (\ref{Q12}), we can obtain that
$$
      m_*-\frac{s_3^{q}}{q}\langle I'(u_0),u_0^+\rangle
\leq  \frac{t_3^{q}}{q}\langle I'(u_0),u_0^-\rangle+m_*
\leq  m_*,
$$
which implies that $\langle I'(u_0),u_0^+ \rangle \geq0$. Similarly, we can also obtain that $\langle I'(u_0),u_0^- \rangle \geq0$. Then by $(\ref{Q11})$ and $(\ref{Q12})$, we have $\langle I'(u_0),u_0^\pm \rangle =0$ and then $\langle I'(u_0),u_0 \rangle =0$. Furthermore, according to $(\ref{Q14})$, we can obtain that $I(u_0)=m_*$ and $u_0\in\mathcal{M}$. Note that $u_0^+\neq0$. If we let $s_3=0$ and $t_3=0$ in (\ref{Eq5}), then we have
$$
m_*=I(u_0)\geq\left(\frac{1}{p}-\frac{1}{q}\right)\|u_0^+\|^p+\left(\frac{1}{p}-\frac{1}{q}\right)\|u_0^-\|^p>0.
$$
 Through arguments similar to the above, we can also conclude that $c_*>0$ can also be achieved.\qed

\vskip2mm
\noindent
{\bf Lemma 2.10.} {\it Assume that $(C_1)$ and $(C_2)$ hold. If $u_0\in\mathcal{M}$ and $I(u_0)=m_*$, then $u_0$ is a critical point of $I$.}
\vskip0mm
\noindent
{\bf Proof.} Arguing by contradiction. If we suppose that $I'(u_0)\neq0$ for all $u_0\in\mathcal{D}$, then there are two positive constants $\delta$ and $\vartheta$ such that
$$
\|I'(u)\|\geq \vartheta,\;\; \forall\;\|u-u_0\|\leq 3\delta.
$$
Since $u_0\in\mathcal{M}$, we have $\langle I'(u_0),u_0^\pm\rangle =0$, and by Lemma 2.2, for all $s,t\geq0$, there exists
\begin{eqnarray}
\label{Q15}
       I(su_0^++tu_0^-)
&\leq& I(u_0)-\sum_{n\in\mathbb{Z}}a(n)\Big|\sum_{i=1}^{\frac{p}{2}-1}\sum_{j=1}^{i-1}2^{i-j}(\Delta u_0^+(n))^{p-(i+j)}
                (\Delta u_0^-(n))^{i+j}\Big|\Theta\nonumber\\
&    & -\left(\frac{1-s^p}{p}-\frac{1-s^{q}}{q}\right)\|u_0^+\|^p
          -\left(\frac{1-t^p}{p}-\frac{1-t^{q}}{q}\right)\|u_0^-\|^p\nonumber\\
& =  & m_*-\sum_{n\in\mathbb{Z}}a(n)\Big|\sum_{i=1}^{\frac{p}{2}-1}\sum_{j=1}^{i-1}2^{i-j}(\Delta u_0^+(n))^{p-(i+j)}
                (\Delta u_0^-(n))^{i+j}\Big|\Theta\nonumber\\
&    & -\left(\frac{1-s^p}{p}-\frac{1-s^{q}}{q}\right)\|u_0^+\|^p
          -\left(\frac{1-t^p}{p}-\frac{1-t^{q}}{q}\right)\|u_0^-\|^p.
\end{eqnarray}
Let $D=(\frac{1}{2},\frac{3}{2})\times(\frac{1}{2},\frac{3}{2})$. It follows from $(\ref{Q15})$ that
\begin{eqnarray}
\label{Q16}
\mathcal{Y}:=\max_{(s,t)\in\partial D}I(su_0^++tu_0^-)<m_*.
\end{eqnarray}
For $\epsilon:=\min\{\frac{m_*-\mathcal{Y}}{3},\frac{\vartheta\delta}{8}\}$ and $S_\delta:=B(u_0,\delta)$, by \cite{Willem 1996}, we can obtain a deformation $\eta\in\mathcal{C}([0,1]\times \mathcal{D},\mathcal{D})$ such that\\
$(\text{i})\;\eta(1,u)=u\;\;\text{if}\;\;|I(u)-m_*|>2\epsilon$;\\
$(\text{ii})\;\eta(1,I^{m_*+\epsilon}\cap S_\delta)\subset I^{m_*-\epsilon}$,\;\;\text{where}\;\; $I^{c}:=\{u\in \mathcal{D}:I(u)\leq c\}$;\\
$(\text{iii})\;I(\eta(1,u))\leq I(u),\;\;\forall u\in{\mathcal{D}}$;\\
$(\text{iv})\;\eta(1,u)\;\;\text{is a homeomorphism of}\;\;\mathcal{D}$.\\
For one thing, in virtue of $(\ref{Q15})$, (iii) and for all $s,t\geq0$ which make $|s-1|^2+|t-1|^2\geq\delta^2/{\|u_0\|^2}$ holds, there exists
\begin{eqnarray}
\label{Q17}
     I(\eta(1,su_0^++tu_0^-))
\leq I(su_0^++tu_0^-)
  <  I(u_0)=m_*.
\end{eqnarray}
For another, by Corollary 2.4, for all $s,t\geq0$, we have that $I(su_0^++tu_0^-)\leq I(u_0)=m_*$. Then by (ii), we have
\begin{eqnarray}
\label{Q18}
I(\eta(1,su_0^++tu_0^-))\leq m_*-\epsilon,\;\;\text{for}\;s,t\geq0,\;|s-1|^2+|t-1|^2<\delta^2/{\|u_0\|^2}.
\end{eqnarray}
In virtue of (\ref{Q16}), (\ref{Q17}) and (\ref{Q18}), we have
\begin{eqnarray}
\label{Q19}
\max_{(s,t)\in\bar{D}}I(\eta(1,su_0^++tu_0^-))<m_*.
\end{eqnarray}
Define $k(s,t)=su_0^++tu_0^-$. Next, we prove that $\eta(1,k(D))\cap\mathcal{M}\neq\varnothing$. Set $\gamma(s,t):=\eta(1,k(s,t))$,
\begin{eqnarray*}
      \phi_1(s,t)
& := & \Big(\left\langle I'(k(s,t)),u_0^+\right\rangle,\left\langle I'(k(s,t)),u_0^-\right\rangle\Big)\\
& := & \big(y_1(s,t),y_2(s,t)\big),
\end{eqnarray*}
and
\begin{eqnarray*}
      \phi_2(s,t)
  :=  \left(\frac{1}{s}\left\langle I'(\gamma(s,t)),(\gamma(s,t))^+\right\rangle,
            \frac{1}{t}\left\langle I'(\gamma(s,t)),(\gamma(s,t))^-\right\rangle\right).
\end{eqnarray*}
Note that $\phi_1(s,t)$ and $\phi_2(s,t)$ are two-dimensional vectors. According to (\ref{Eq16}) and (\ref{Eq17}), it is obvious that $y_1(s,t)=\frac{1}{s}h_1(s,t)$ and $y_2(s,t)=\frac{1}{t}h_2(s,t)$. So $\phi_1(s,t)$ is a $C^1$ function of $s,t$ and we have
\begin{eqnarray*}
      \frac{\partial y_1(s,t)}{\partial s}|_{(1,1)}
& = & (p-1)\|u_0^+\|^p-\sum_{n\in\mathbb{Z}}(q-1)c(n)|u_0^+(n)|^q\ln{|u_0^+(n)|^r}-\sum_{n\in\mathbb{Z}}rc(n)|u_0^+(n)|^q\\
&   & +(p-i-j-1)\sum_{n\in\mathbb{Z}}a(n)\Big|\sum_{i=1}^{\frac{p}{2}-1}\sum_{j=0}^i C_{\frac{p}{2}-1}^i C_i^j2^{i-j}
                 (\Delta u_0^+(n))^{p-(i+j)}(\Delta u_0^-(n))^{i+j}\Big|\\
&   &  +(p-i-j-1)\sum_{n\in\mathbb{Z}}a(n)\Big|\sum_{i=1}^{\frac{p}{2}}\sum_{j=0}^{i-1} C_{\frac{p}{2}-1}^{i-1}
                  C_{i-1}^{j}2^{i-1-j}(\Delta u_0^+(n))^{p-(i+j)}(\Delta u_0^-(n))^{i+j}\Big|,
\end{eqnarray*}
\begin{eqnarray*}
      \frac{\partial y_1(s,t)}{\partial t}|_{(1,1)}
& = & (i+j)\sum_{n\in\mathbb{Z}}a(n)\Big|\sum_{i=1}^{\frac{p}{2}-1}\sum_{j=0}^i C_{\frac{p}{2}-1}^i C_i^j2^{i-j}
                 (\Delta u_0^+(n))^{p-(i+j)}(\Delta u_0^-(n))^{i+j}\Big|\\
&   &  +(i+j)\sum_{n\in\mathbb{Z}}a(n)\Big|\sum_{i=1}^{\frac{p}{2}}\sum_{j=0}^{i-1} C_{\frac{p}{2}-1}^{i-1}
                 C_{i-1}^{j}2^{i-1-j}(\Delta u_0^+(n))^{p-(i+j)}(\Delta u_0^-(n))^{i+j}\Big|.
\end{eqnarray*}
Similarly, we also have
\begin{eqnarray*}
       \frac{\partial y_2(s,t)}{\partial s}|_{(1,1)}
& = &  (p-i-j)\sum_{n\in\mathbb{Z}}a(n)\Big|\sum_{i=1}^{\frac{p}{2}-1}\sum_{j=1}^{i} C_{\frac{p}{2}-1}^{i-1}
               C_{i-1}^{j-1}2^{i-j}(\Delta u_0^+(n))^{p-(i+j)}(\Delta u_0^-(n))^{i+j}\Big|\\
&   &   +(\frac{p}{2}-j)\sum_{n\in\mathbb{Z}}a(n)\Big|\sum_{j=1}^{\frac{p}{2}-1}C_{\frac{p}{2}-1}^{j-1}2^{\frac{p}{2}-j}
               (\Delta u_0^+(n))^{\frac{p}{2}-j}(\Delta u_0^-(n))^{\frac{p}{2}+j}\Big|\\
&   &   +(p-i-j)\sum_{n\in\mathbb{Z}}a(n)\Big|\sum_{i=1}^{\frac{p}{2}}\sum_{j=0}^{i-1} C_{\frac{p}{2}-1}^{i-1}
               C_{i-1}^{j}2^{i-1-j}(\Delta u_0^+(n))^{p-(i+j)}(\Delta u_0^-(n))^{i+j}\Big|,
\end{eqnarray*}
\begin{eqnarray*}
      \frac{\partial y_2(s,t)}{\partial t}|_{(1,1)}
& = & (p-1)\|u_0^-\|^p-\sum_{n\in\mathbb{Z}}(q-1)c(n)|u_0^-(n)|^q\ln{|u_0^-(n)|^r}-\sum_{n\in\mathbb{Z}}rc(n)|u_0^-(n)|^q\\
&   & +(i+j-1)\sum_{n\in\mathbb{Z}}a(n)\Big|\sum_{i=1}^{\frac{p}{2}-1}\sum_{j=1}^{i} C_{\frac{p}{2}-1}^{i-1}
               C_{i-1}^{j-1}2^{i-j}(\Delta u_0^+(n))^{p-(i+j)}(\Delta u_0^-(n))^{i+j}\Big|\\
&   &   +(\frac{p}{2}+j-1)\sum_{n\in\mathbb{Z}}a(n)\Big|\sum_{j=1}^{\frac{p}{2}-1}C_{\frac{p}{2}-1}^{j-1}2^{\frac{p}{2}-j}
               (\Delta u_0^+(n))^{\frac{p}{2}-j}(\Delta u_0^-(n))^{\frac{p}{2}+j}\Big|\\
&   &   +(i+j-1)\sum_{n\in\mathbb{Z}}a(n)\Big|\sum_{i=1}^{\frac{p}{2}}\sum_{j=0}^{i-1} C_{\frac{p}{2}-1}^{i-1}
               C_{i-1}^{j}2^{i-1-j}(\Delta u_0^+(n))^{p-(i+j)}(\Delta u_0^-(n))^{i+j}\Big|.
\end{eqnarray*}
Let
\begin{gather*}
M=\begin{bmatrix}
  \frac{\partial y_1(s,t)}{\partial s}|_{(1,1)} & \frac{\partial y_2(s,t)}{\partial s}|_{(1,1)} \\
  \frac{\partial y_1(s,t)}{\partial t}|_{(1,1)} & \frac{\partial y_2(s,t)}{\partial t}|_{(1,1)}
  \end{bmatrix}.\quad
  \end{gather*}
By $(\ref{Q9})$ and $(\ref{Q10})$, we have
\begin{small}
\begin{eqnarray*}
        \frac{\partial y_1(s,t)}{\partial s}|_{(1,1)}
& = & (p-1)\|u_0^+\|^p-\sum_{n\in\mathbb{Z}}rc(n)|u_0^+(n)|^q \\
&   & +(p-i-j-1)\sum_{n\in\mathbb{Z}}a(n)\Big|\sum_{i=1}^{\frac{p}{2}-1}\sum_{j=0}^i C_{\frac{p}{2}-1}^i C_i^j2^{i-j}
                 (\Delta u_0^+(n))^{p-(i+j)}(\Delta u_0^-(n))^{i+j}\Big|\\
&   &  +(p-i-j-1)\sum_{n\in\mathbb{Z}}a(n)\Big|\sum_{i=1}^{\frac{p}{2}}\sum_{j=0}^{i-1} C_{\frac{p}{2}-1}^{i-1}
                  C_{i-1}^{j}2^{i-1-j}(\Delta u_0^+(n))^{p-(i+j)}(\Delta u_0^-(n))^{i+j}\Big|\\
&   & -(q-1)\|u_0^+\|^p-(q-1)\sum_{n\in\mathbb{Z}}a(n)\Big|\sum_{i=1}^{\frac{p}{2}-1}\sum_{j=0}^i C_{\frac{p}{2}-1}^i C_i^j2^{i-j}(\Delta u_0^+(n))^{p-(i+j)}(\Delta u_0^-(n))^{i+j}\Big|\\
&   &   -(q-1)\sum_{n\in\mathbb{Z}}a(n)\Big|\sum_{i=1}^{\frac{p}{2}}\sum_{j=0}^{i-1} C_{\frac{p}{2}-1}^{i-1}
                  C_{i-1}^{j}2^{i-1-j}(\Delta u_0^+(n))^{p-(i+j)}(\Delta u_0^-(n))^{i+j}\Big|\\
& = & (p-q)\left[\|u_0^+\|^p
        +\sum_{n\in\mathbb{Z}}a(n)\Big|\sum_{i=1}^{\frac{p}{2}-1}\sum_{j=0}^i C_{\frac{p}{2}-1}^i C_i^j2^{i-j}
                 (\Delta u_0^+(n))^{p-(i+j)}(\Delta u_0^-(n))^{i+j}\Big|\right. \\
&   &  \left.+\sum_{n\in\mathbb{Z}}a(n)\Big|\sum_{i=1}^{\frac{p}{2}}\sum_{j=0}^{i-1} C_{\frac{p}{2}-1}^{i-1}
                  C_{i-1}^{j}2^{i-1-j}(\Delta u_0^+(n))^{p-(i+j)}(\Delta u_0^-(n))^{i+j}\Big|\right]
                  -\sum_{n\in\mathbb{Z}}rc(n)|u_0^+(n)|^q\\
&   & +(-i-j)\sum_{n\in\mathbb{Z}}a(n)\Big|\sum_{i=1}^{\frac{p}{2}-1}\sum_{j=0}^i C_{\frac{p}{2}-1}^i C_i^j2^{i-j}
                 (\Delta u_0^+(n))^{p-(i+j)}(\Delta u_0^-(n))^{i+j}\Big|\\
&   & +(-i-j)\sum_{n\in\mathbb{Z}}a(n)\Big|\sum_{i=1}^{\frac{p}{2}}\sum_{j=0}^{i-1} C_{\frac{p}{2}-1}^{i-1}
                  C_{i-1}^{j}2^{i-1-j}(\Delta u_0^+(n))^{p-(i+j)}(\Delta u_0^-(n))^{i+j}\Big|\\
& = & -\frac{\partial y_1(s,t)}{\partial t}|_{(1,1)}+C^+,
\end{eqnarray*}
\end{small}
and
\begin{eqnarray*}
     \frac{\partial y_2(s,t)}{\partial t}|_{(1,1)}= -\frac{\partial y_2(s,t)}{\partial s}|_{(1,1)}+C^-,
\end{eqnarray*}
where
\begin{eqnarray*}
       C^{+}
& = & (p-q)\left[\|u_0^+\|^p
        +\sum_{n\in\mathbb{Z}}a(n)\Big|\sum_{i=1}^{\frac{p}{2}-1}\sum_{j=0}^i C_{\frac{p}{2}-1}^i C_i^j2^{i-j}
                 (\Delta u_0^+(n))^{p-(i+j)}(\Delta u_0^-(n))^{i+j}\Big|\right. \\
&   &  \left.+\sum_{n\in\mathbb{Z}}a(n)\Big|\sum_{i=1}^{\frac{p}{2}}\sum_{j=0}^{i-1} C_{\frac{p}{2}-1}^{i-1}
                  C_{i-1}^{j}2^{i-1-j}(\Delta u_0^+(n))^{p-(i+j)}(\Delta u_0^-(n))^{i+j}\Big|\right]
                  -\sum_{n\in\mathbb{Z}}rc(n)|u_0^+(n)|^q
\end{eqnarray*}
and
\begin{eqnarray*}
       C^-
& = & (p-q)\left[\|u_0^-\|^p
         +\sum_{n\in\mathbb{Z}}a(n)\Big|\sum_{i=1}^{\frac{p}{2}-1}\sum_{j=1}^{i} C_{\frac{p}{2}-1}^{i-1}
               C_{i-1}^{j-1}2^{i-j}(\Delta u_0^+(n))^{p-(i+j)}(\Delta u_0^-(n))^{i+j}\Big|\right.\\
&   &   \left.+\sum_{n\in\mathbb{Z}}a(n)\Big|\sum_{j=1}^{\frac{p}{2}-1}C_{\frac{p}{2}-1}^{j-1}2^{\frac{p}{2}-j}
               (\Delta u_0^+(n))^{\frac{p}{2}-j}(\Delta u_0^-(n))^{\frac{p}{2}+j}\Big|\right.\\
&   &   \left.+\sum_{n\in\mathbb{Z}}a(n)\Big|\sum_{i=1}^{\frac{p}{2}}\sum_{j=0}^{i-1} C_{\frac{p}{2}-1}^{i-1}
               C_{i-1}^{j}2^{i-1-j}(\Delta u_0^+(n))^{p-(i+j)}(\Delta u_0^-(n))^{i+j}\Big|\right]
         -\sum_{n\in\mathbb{Z}}rc(n)|u_0^-(n)|^q
\end{eqnarray*}
It follows from $1<p<q$ and $u_0^\pm\neq 0$ that $C^+<0$, $C^-<0$, $\frac{\partial y_1(s,t)}{\partial t}|_{(1,1)}\geq0$ and $\frac{\partial y_2(s,t)}{\partial s}|_{(1,1)}\geq0$. Then, we have
\begin{eqnarray*}
      \det M
& = & \frac{\partial y_1(s,t)}{\partial s}|_{(1,1)}\times\frac{\partial y_2(s,t)}{\partial t}|_{(1,1)}
          -\frac{\partial y_1(s,t)}{\partial t}|_{(1,1)}\times\frac{\partial y_2(s,t)}{\partial s}|_{(1,1)}\\
& = & \left(-\frac{\partial y_1(s,t)}{\partial t}|_{(1,1)}+C^+\right)
            \times\left(-\frac{\partial y_2(s,t)}{\partial s}|_{(1,1)}+C^-\right)
          -\frac{\partial y_1(s,t)}{\partial t}|_{(1,1)}\times\frac{\partial y_2(s,t)}{\partial s}|_{(1,1)}\\
& = & -C^-\cdot\frac{\partial y_1(s,t)}{\partial t}|_{(1,1)}
       -C^+\cdot\frac{\partial y_2(s,t)}{\partial s}|_{(1,1)}+C^+C^-\\
& > & 0,
\end{eqnarray*}
which implies that $\det M\neq 0$. According to the topological degree theory \cite{Chang 2005}, we can obtain that $\text{deg}(\phi_1,D,(0,0))=1$. In virtue of (\ref{Q16}) and (i), there exists $\gamma=k$ on $\partial D$. As a consequence, it follows from homotopy type invariance of the topological degree theory that
$$
\text{deg}(\phi_2,D,(0,0))=\text{deg}(\phi_1,D,(0,0))=1,
$$
which implies that $\phi_2(s_4,t_4)=0$ for some $(s_4,t_4)\in D$ and so $\eta(1,k(s_4,t_4))=\gamma(s_4,t_4)\in\mathcal{M}$. Then by the definition of $\mathcal{M}$, we know that $I(\eta(1,k(s_4,t_4)))\ge m^*$. This contradicts with (\ref{Q19}). Hence, $I'(u_0)=0$, that is, $u_0$ is a critical point of $I$.
\qed
\vskip0mm
\noindent
\vskip2mm
{\section{The existence of sign-changing solutions}}
\setcounter{equation}{0}
\vskip0mm
\noindent
In this section, we will prove the existence of sign-changing solutions which only changes sign one times.
\vskip0mm
\noindent
{\bf Proof of Theorem 1.1:} Firstly, it follows from Lemma 2.9 and Lemma 2.10 that problem $(\ref{eq1})$ has  a sign-changing solution $u_0\in\mathcal{M}$ such that
\begin{eqnarray}
\label{W2}
I(u_0)=m_*\;\;\text{and}\;\;I'(u_0)=0.
\end{eqnarray}
Next we prove that $u_0$ only changes sign one times. Denote $u_0=u_1+u_2+u_3$, where
\begin{eqnarray}
\label{W3}
u_1\geq0,\;\;u_2\leq0,\;\;V_1\cap V_2=\emptyset,\;\;
u_1\mid_{\mathbb{Z}\setminus(V_1\cup V_2)}=u_2\mid_{\mathbb{Z}\setminus(V_1\cup V_2)}=u_3\mid_{V_1\cup V_2}=0,
\end{eqnarray}
$$
V_1:=\{n\in\mathbb{Z}:u_1(n)>0\},\;\;V_2:=\{n\in\mathbb{Z}:u_2(n)<0\},
$$
and $V_1=\{n_1,n_1+1,\cdots,n_1+m_1\}$, $V_2=\{n_2,n_2+1,\cdots,n_2+m_2\}$, where the value of $n_1$ or $n_2$ may be $-\infty$ and the value of $n_1+m_1$ or $n_2+m_2$ may be $+\infty$.
\par
Setting $w=u_1+u_2$, it is easy to see that $w^+=u_1$, $w^-=u_2$ and $w^\pm\neq0$. According to Lemma 2.7, there is a unique pair of positive constants $s_4,t_4$ such that $s_4w^++t_4w^-\in\mathcal{M}$. In virtue of $I'(u_0)=0$, we can derive that $\langle I'(u_0),w^\pm\rangle=0$. Then by $(\ref{eq10})$, we can obtain that
\begin{eqnarray}
\label{W33}
      \langle I'(u_0),w^+\rangle
& = & \sum_{n\in\mathbb{Z}}\left[a(n)|\Delta w(n)+\Delta u_3(n)|^{p-2}(\Delta w(n)+\Delta u_3(n))\Delta w^+(n)\right]\nonumber\\
&   & +\sum_{n\in\mathbb{Z}}\left[b(n)|w(n)+u_3(n)|^{p-2}(w(n)+u_3(n))w^+(n)\right] \nonumber\\
&   & -\sum_{n\in\mathbb{Z}}c(n)|w(n)+u_3(n)|^{q-2}(w(n)+u_3(n))w^+(n)\ln{|w(n)+u_3(n)|^r}\nonumber\\
& = & \langle I'(w),w^+ \rangle+\sum_{n\in V_1}a(n)|\Delta w(n)+\Delta u_3(n)|^{p-2}\Delta u_3(n)\Delta w^+(n)\nonumber\\
&   & +\sum_{n\in {\mathbb{Z}\setminus V_1}}a(n)|\Delta w(n)+\Delta u_3(n)|^{p-2}\Delta u_3(n)\Delta w^+(n)\nonumber\\
&   & +\sum_{n\in\mathbb{Z}}a(n)\Big|\sum_{i=1}^{\frac{p}{2}-1} C_{\frac{p}{2}-1}^i
                 (\Delta w(n))^{p-2-i}(\Delta u_3(n))^{i}\Big|\Delta w(n)\Delta w^+(n)\nonumber\\
& = & \langle I'(w),w^+ \rangle-a(n_1-1)\Big|w^+(n_1)-u_3(n_1-1)\Big|^{p-2}w^+(n_1)u_3(n_1-1)\nonumber\\
&   & -a(n_1+m_1)\Big|-w^+(n_1+m_1)+u_3(n_1+m_1+1)\Big|^{p-2}w^+(n_1+m_1)u_3(n_1+m_1+1)\nonumber\\
&   & +\sum_{n\in\mathbb{Z}}a(n)\Big|\sum_{i=1}^{\frac{p}{2}-1} C_{\frac{p}{2}-1}^i
                 (\Delta w(n))^{p-2-i}(\Delta u_3(n))^{i}\Big|\Delta w(n)\Delta w^+(n).
\end{eqnarray}
Note that
\begin{eqnarray*}
     \Delta w(n)\Delta w^+(n)
&=&  \Delta w^+(n)\Delta w^+(n)+\Delta w^-(n)\Delta w^+(n)\\
&=&  (\Delta w^+(n))^2+\Delta w^-(n)\Delta w^+(n)\\
&=&  (\Delta w^+(n))^2-w^-(n+1)w^+(n)-w^-(n)w^+(n+1)\geq0.
\end{eqnarray*}
According to (\ref{W33}), one has
\begin{eqnarray}
\label{W4}
       \langle I'(w),w^+\rangle
& = & a(n_1-1)\Big|w^+(n_1)-u_3(n_1-1)\Big|^{p-2}w^+(n_1)u_3(n_1-1)\nonumber\\
&   & +a(n_1+m_1)\Big|-w^+(n_1+m_1)+u_3(n_1+m_1+1)\Big|^{p-2}w^+(n_1+m_1)u_3(n_1+m_1+1)\nonumber\\
&   & -\sum_{n\in\mathbb{Z}}a(n)\Big|\sum_{i=1}^{\frac{p}{2}-1} C_{\frac{p}{2}-1}^i
                 (\Delta w(n))^{p-2-i}(\Delta u_3(n))^{i}\Big|\Delta w(n)\Delta w^+(n)\nonumber\\
&\leq& 0.
\end{eqnarray}
Similarly, we can obtain that
\begin{eqnarray}
\label{W44}
       \langle I'(w),w^-\rangle
& = & a(n_2-1)\Big|w^-(n_2)-u_3(n_2-1)\Big|^{p-2}w^-(n_2)u_3(n_2-1)\nonumber\\
&   & +a(n_2+m_2)\Big|-w^-(n_2+m_2)+u_3(n_2+m_2+1)\Big|^{p-2}w^-(n_2+m_2)u_3(n_2+m_2+1)\nonumber\\
&   & -\sum_{n\in\mathbb{Z}}a(n)\Big|\sum_{i=1}^{\frac{p}{2}-1} C_{\frac{p}{2}-1}^i
                 (\Delta w(n))^{p-2-i}(\Delta u_3(n))^{i}\Big|\Delta w(n)\Delta w^-(n)\nonumber\\
&\leq& 0.
\end{eqnarray}
On the basis of $(\ref{eq7})$, $(\ref{eq10})$, $(\ref{Eq5})$, $(\ref{W2})$, $(\ref{W3})$, $(\ref{W4})$ and $(\ref{W44})$, using the same processing method as (2.7), we have
\begin{small}
\begin{eqnarray*}
      m_*
& = & I(u_0)-\frac{1}{q}\langle I'(u_0),u_0\rangle\\
& = & I(w)+I(u_3)+\left(\frac{1}{p}-\frac{1}{q}\right)\sum_{n\in\mathbb{Z}}a(n)\Big|\sum_{i=1}^{\frac{p}{2}-1}\sum_{j=0}^i
             C_{\frac{p}{2}}^iC_i^j2^{i-j}(\Delta w(n))^{p-(i+j)}(\Delta u_3(n))^{i+j}\Big|\\
&   & -\frac{1}{q}\langle I'(w),w\rangle-\frac{1}{q}\langle I'(u_3),u_3\rangle+\left(\frac{1}{p}-\frac{1}{q}\right)
             \sum_{n\in\mathbb{Z}}a(n)\Big|\sum_{j=0}^{\frac{p}{2}-1}C_{\frac{p}{2}}^j2^{\frac{p}{2}-j}
             (\Delta w(n))^{\frac{p}{2}-j}(\Delta u_3(n))^{\frac{p}{2}+j}\Big|\\
&\geq& I(s_0w^++t_0w^-)+\frac{1-s_0^q}{q}\langle I'(w),w^+\rangle+\frac{1-t_0^q}{q}\langle I'(w),w^-\rangle
       +I(u_3)-\frac{1}{q}\langle I'(w),w\rangle-\frac{1}{q}\langle I'(u_3),u_3\rangle\\
& = & I(s_0w^++t_0w^-)-\frac{s_0^q}{q}\langle I'(w),w^+\rangle-\frac{t_0^q}{q}\langle I'(w),w^-\rangle
       +I(u_3)-\frac{1}{q}\langle I'(u_3),u_3\rangle\\
& = & I(s_0w^++t_0w^-)-\frac{s_0^q}{q}\langle I'(w),w^+\rangle-\frac{t_0^q}{q}\langle I'(w),w^-\rangle
         +\left(\frac{1}{p}-\frac{1}{q}\right)\|u_3\|+\frac{r}{q^2}\sum_{n\in\mathbb{Z}}c(n)|u_3(n)|^q\\
&\geq& m_*+\left(\frac{1}{p}-\frac{1}{q}\right)\|u_3\|,
\end{eqnarray*}
\end{small}
which implies that $u_3=0$. Thus, $u_0$ only changes the sign one times.
\qed
\vskip2mm
{\section{The existence of ground state solutions}}
\setcounter{equation}{0}
\vskip0mm
\noindent
In this section, we will prove the existence of Nehari type ground state solution for $(\ref{eq1})$ and provide the relationship between sign-changing ground state energy and the ground state energy. We mainly use the method in \cite{Chen 2013,Chang 2023} to prove that the functional $I$ satisfies the Cerami condition at any level $d\in(0,\infty)$, and then use the method in  \cite{Chen 2018} to prove that the functional $I$ has a mountain pass geometry. To prove the above conclusions, we need the following lemmas.
\vskip2mm
\noindent
{\bf Lemma 4.1.} \cite{Schechter 1991} {\it Let $X$ be a real Banach space. For some constants $\alpha,\beta,\rho>0$ and $e\in X$ with $\|e\|_X>\rho$, there exists a functional $I\in C^1(X,\mathbb{R})$ satisfying the following mountain pass geometry:
$$
\max\big\{I(0),I(e)\big\}\leq\alpha<\beta\leq\inf_{\|u\|_X=\rho} I(u).
$$
Set $d_0=\inf_{\gamma\in\Gamma}\max_{0\leq t\leq 1} I\big(\gamma(t)\big)$, where $\Gamma=\big\{\gamma\in C([0,1],X):\gamma(0)=0\;\;\text{and}\;\;\gamma(1)=e\big\}$. Then there exists a Cerami sequence $\{u_k\}\subset X$ of $I$ at level $d_0$, where  a sequence $\{u_k\}$ is called   Cerami sequence at a level $d_0$ if it satisfies}
\begin{eqnarray}
\label{pp2}
I(u_k)\rightarrow d_0\;\;\text{\it and}\;\; \|I'(u_k)\|(1+\|u_k\|)\rightarrow 0.
\end{eqnarray}
{\bf Remark 4.1.} It is easy to obtain that $d_0\ge \beta>0$ (for example, see the proof of Theorem 1.15 in \cite{Willem 1996}).

\vskip2mm
\noindent
{\bf Lemma 4.2.} {\it The Cerami sequence $\{u_k\}\subset  \mathcal{D}$ at any level $d_0\in (0,+\infty)$
has at least one convergent subsequence in} $\mathcal{D}$.\\
{\bf Proof.} Since   $\{u_k\}$ is a Cerami sequence at the level $d_0$, then (\ref{pp2}) hold. We claim that $\{u_k\}$ is bounded in $\mathcal{D}$. Arguing by contradiction, we suppose that $\{u_k\}$ is not bounded in $\mathcal{D}$, that is, there appears a subsequence, still denoted by $\{u_k\}$, such that $\|u_k\|\rightarrow+\infty$ as $k\rightarrow+\infty$.
\par
Let $w_k=\frac{u_k}{\|u_k\|}$. Then there exists a subsequence, still denoted by $\{w_k\}$, and a function $w\in \mathcal{D}$ such that
\begin{eqnarray}
\label{W5}
\begin{cases}
 w_k   \rightharpoonup w    &\text{in}\;\mathcal{D},\\
 w_k(n)\rightarrow     w(n) &\mbox{for each } n\in \mathbb{Z},\\
 w_k   \rightarrow     w    &\text{in} \; l^\kappa(\mathbb{Z},\mathbb{R}), \kappa\in[p,+\infty].
 \end{cases}
\end{eqnarray}
Then we will prove the claim by discussing the following two cases.
\par
{\bf Case 1:} $w=0$.
\par
Set $t_k\in[0,1]$ such that $I(t_ku_k)=\max_{t\in[0,1]} I(tu_k)$. For any given constants $\tau>0$ and $N>0$, it follows from the unboundedness of $\{u_k\}$ that
$$
\|u_k\|> (p\tau+1)^{\frac{1}{p}},\;\;\text{for large enough}\;k\geq N.
$$
Set $\bar{w}_k=(p\tau+1)^{\frac{1}{p}}w_k$. In virtue of $(C_2)$, the boundedness of $\{w_k\}$, (\ref{Eq1}), $(\ref{eq9})$ and Lebesgue dominate convergence theorem, we have
\begin{eqnarray*}
&   & \lim_{k\rightarrow+\infty}\sum_{n\in\mathbb{Z}}c(n)|\bar{w}_k(n)|^q\ln{|\bar{w}_k(n)|^r}\\
& = & \lim_{k\rightarrow+\infty}\sum_{n\in\mathbb{Z}}c(n)(p\tau+1)^{\frac{q}{p}}|w_k(n)|^q\ln{((p\tau+1)^\frac{r}{p}|w_k(n)|^r})\\
& = & \lim_{k\rightarrow +\infty} (p\tau+1)^{\frac{q}{p}}\left(\sum_{n\in\mathbb{Z}}c(n)|w_k(n)|^q\ln {(p\tau+1)^\frac{r}{p}}
        +\sum_{n\in\mathbb{Z}}c(n)|w_k(n)|^q\ln {|w_k(n)|^r}\right)\\
& = & (p\tau+1)^{\frac{q}{p}}\left(\sum_{n\in\mathbb{Z}}c(n)|w(n)|^q\ln {(p\tau+1)^\frac{r}{p}}
        +\sum_{n\in\mathbb{Z}}c(n)|w(n)|^q\ln {|w(n)|^r}\right)\\
& = & 0.
\end{eqnarray*}
Then for $k$ large enough, we can obtain that
\begin{eqnarray*}
       I(t_k u_k)\geq I\left(\frac{(p\tau+1)^{\frac{1}{p}}}{\|u_k\|}u_k\right)=I(\bar{w}_k)
\geq \frac{1}{p}\|\bar{w}_k\|^p-\frac{1}{q}\sum_{n\in\mathbb{Z}}c(n)|\bar{w}_k(n)|^q\ln {|\bar{w}_k(n)|^r}\geq\tau.
\end{eqnarray*}
According to the arbitrariness of $\tau$, we can obtain that
\begin{eqnarray}
\label{W6}
\lim_{k\rightarrow +\infty} I(t_k u_k)=+\infty.
\end{eqnarray}
If $t_k=1$, substituting it into (\ref{W6}) can obtain $\lim_{k\rightarrow \infty} I(u_k)=+\infty$, which contradicts with (\ref{pp2}). And then it follows from $I(0)=0$ that $t_k\in(0,1)$. Thus $\frac{d}{dt} I(tu_k)\mid_{t=t_k}=0$. Therefore, according to the definition of $\{u_k\}$, we can obtain that
\begin{eqnarray*}
       I(t_k u_k)
& =  & I(t_k u_k)-\frac{1}{q}\langle I'(t_k u_k),t_k u_k\rangle
           =\left(\frac{1}{p}-\frac{1}{q}\right)\|t_k u_k\|^p+\frac{r}{q^2}\sum_{n\in\mathbb{Z}}c(n)|t_k u_k(n)|^q\\
&\leq& \left(\frac{1}{p}-\frac{1}{q}\right)\|u_k\|^p+\frac{r}{q^2}\sum_{n\in\mathbb{Z}}c(n)|u_k(n)|^q
           =I(u_k)-\frac{1}{q}\langle I'(u_k),u_k\rangle\\
&\leq& d+o(1),
\end{eqnarray*}
which is contrary to $(\ref{W6})$. So the assumption is not valid, that is, $\{u_k\}$ is bounded in $\mathcal{D}$.
\par
{\bf Case 2:} $w\neq 0$.
\par
Let $V'=\big\{n\in\mathbb{Z};\;w\neq 0\big\}$. Then $|u_k(n)|\to +\infty$ as $k\rightarrow+\infty$ for each $n\in V'$. According to the fact $\|u_k\|\to +\infty$, as $k\to+\infty$, and $I(u_k)\leq c_*$, there holds $\frac{I(u_k)}{\|u_k\|^p}\to 0$, as $k\to+\infty$, that is
$$
  \frac{1}{p}+\frac{r\sum_{n\in\mathbb{Z}}c(n)|u_k(n)|^q}{q^2\|u_k\|^p}
    -\frac{\sum_{n\in\mathbb{Z}}c(n)|u_k(n)|^q\ln{|u_k(n)|^r}}{q\|u_k\|^p}
=  o_k(1),
$$
which together with the definition of $\mathcal{D}$, $(C_1)$ and $(C_2)$ implies that $G(u_k):=\frac{r\sum_{n\in\mathbb{Z}}c(n)|u_k(n)|^q}{q^2\|u_k\|^p}
-\frac{\sum_{n\in\mathbb{Z}}c(n)|u_k(n)|^q\ln{|u_k(n)|^r}}{q\|u_k\|^p}$ is bounded, that is, there is a positive constant $M_2$ such that \begin{eqnarray}\label{aaa1}
|G(u_k)|\le M_2, \ \forall \ k\in\mathbb N.
\end{eqnarray}  We set
\begin{eqnarray}
\label{W7}
G(u_k)=\sum_{n\in \mathbb{Z}\setminus V';\;|u_k(n)|\leq M_3}(\cdot)
          +\sum_{n\in \mathbb{Z}\setminus V';\;|u_k(n)|>M_3}(\cdot)+\sum_{n\in V'}(\cdot)
   :=I+II+III,
\end{eqnarray}
where  $M_3=e^{\frac{1}{q}}>0$ and $(\cdot)=\frac{rc(n)|u_k(n)|^q}{q^2\|u_k\|^p}
-\frac{c(n)|u_k(n)|^q\ln{|u_k(n)|^r}}{q\|u_k\|^p}$. For $I$ in (\ref{W7}), according to $(\ref{eq9})$, Lemma 2.1, (\ref{Q111}) and $\sum_{n\in\mathbb{Z}}c(n)<+\infty$, there are three positive constants $\varepsilon_6$ and $C_{\varepsilon_6}$ such that
\begin{eqnarray}
\label{W8}
       I
&  = & \sum_{n\in\mathbb{Z}\setminus V';\;|u_k(n)|\leq M_3}
             \left(\frac{rc(n)|u_{k}(n)|^q}{q^2\|u_{k}\|^p}
                   -\frac{c(n)|u_{k}(n)|^q\ln{|u_{k}(n)|^r}}{q\|u_{k}\|^p}\right)\nonumber\\
&\leq& \sum_{n\in\mathbb{Z}\setminus V';\;|u_{k}(n)|\leq M_3}
             \left(\frac{rc(n)|u_{k}(n)|^q}{q^2\|u_{k}\|^p}
                   +\frac{c(n)\varepsilon_6|u_{k}(n)|^p+c(n)C_{\varepsilon_6}|u_{k}(n)|^\zeta}{q\|u_{k}\|^p}\right)\nonumber\\
&\leq& \sum_{n\in\mathbb{Z}\setminus V';\;|u_{k}(n)|\leq M_3}
             \left(\frac{rc(n)M_3}{q^2\|u_{k}\|^p}
                   +\frac{c(n)\varepsilon_6 M_3+c(n)C_{\varepsilon_6}M_3}{q\|u_{k}\|^p}\right)\nonumber\\
&\leq&
             \left(\frac{rM_3\sum_{n\in\mathbb{Z}}c(n)}{q^2\|u_{k}\|^p}
                   +\frac{(\varepsilon_6+C_{\varepsilon_6})M_3\sum_{n\in\mathbb{Z}}c(n)}{q\|u_{k}\|^p}\right)\nonumber\\
& \to  & 0, \ \mbox{as }k\to+\infty.
\end{eqnarray}
For $II$ in (\ref{W7}), we have
\begin{eqnarray}
\label{WW99}
  II  & =  &  \sum_{n\in\mathbb{Z}\setminus V';\;|u_{k}(n)|>M_3}
             \left(\frac{rc(n)|u_{k}(n)|^q}{q^2\|u_{k}\|^p}
                   -\frac{c(n)|u_{k}(n)|^q\ln{|u_{k}(n)|^r}}{q\|u_{k}\|^p}\right)\nonumber\\
& =  & \sum_{n\in\mathbb{Z}\setminus V';\;|u_{k}(n)|>M_3}
                 \frac{c(n)|u_{k}(n)|^q}{q^2\|u_{k}\|^p}\left(r-q\ln{|u_{k}(n)|^r}\right)\nonumber\\
&\leq& \sum_{n\in\mathbb{Z}\setminus V';\;|u_{k}(n)|>M_3}
                 \frac{c(n)|w_k(n)|^p|u_{k}(n)|^{q-p}}{q^2}\left(r-q\ln{M_3^r}\right)\to-\infty, \ \mbox{as }k\to+\infty.
\end{eqnarray}
Note that $|u_k(n)|\to +\infty$ as $k\rightarrow+\infty$ for each $n\in V'$. Then similar to the argument of $II$, we also have
\begin{eqnarray}
\label{WW999}
III= \sum_{n\in V'}\left(\frac{rc(n)|u_{k}(n)|^q}{q^2\|u_{k}\|^p}
           -\frac{c(n)|u_{k}(n)|^q\ln{|u_{k}(n)|^r}}{q\|u_{k}\|^p}\right)\to -\infty, \ \mbox{as }k\to+\infty.
\end{eqnarray}
Thus, $\lim_{k\to\infty}G(u_k)=-\infty$, which contradicts to $(\ref{aaa1})$.
Thus, we deduce that $\{u_k\}$ is bounded in $\mathcal{D}$. Then there exists a subsequence, still denoted by $\{u_k\}$, and a function $u\in \mathcal{D}$ such that
\begin{eqnarray}
\label{W10}
\begin{cases}
 u_k   \rightharpoonup u    &\text{in}\;\mathcal{D},\\
 u_k(n)\rightarrow     u(n) &\mbox{for each } n\in \mathbb{Z},\\
 u_k   \rightarrow     u    &\text{in} \; l^\kappa(\mathbb{Z},\mathbb{R}), \kappa\in[p,+\infty].
 \end{cases}
\end{eqnarray}
Note that $\{u_k\}$ is a Cerami sequence. Then there holds
\begin{eqnarray}
\label{W11}
\lim_{k\to +\infty} \langle I'(u_k),u_k-u\rangle=0.
\end{eqnarray}
Moreover, by $(\ref{W10})$, we have
\begin{eqnarray}
\label{W15}
\lim_{k\rightarrow+\infty} \langle I'(u),u_k-u\rangle =0.
\end{eqnarray}
Note that $\{u_k\}$ is bounded in $E$. It follows from (\ref{eq9}), (\ref{Eq1}),  (\ref{W10}) and $\sum_{n\in\mathbb{Z}}c(n)<+\infty$ that there exist two positive constants $\varepsilon$ and $C_{\varepsilon}$ such that
\begin{eqnarray}
\label{W13}
&    &  \sum_{n\in\mathbb{Z}}c(n)|u(n)|^{q-2}u(n)(u_k(n)-u(n))\ln{|u(n)|^r}\nonumber\\
&\leq&  \sum_{n\in\mathbb{Z}}c(n)(u_k(n)-u(n))
            \left(\varepsilon|u(n)|^{p-1}+C_{\varepsilon}|u(n)|^{\zeta-1}\right)\nonumber\\
& \to  &  0\ \ \mbox{ as } k\to+\infty.
\end{eqnarray}
Similarly, it follows from the boundedness of $\{\|u_k\|\}$, (\ref{eq9}), (\ref{Eq1}),  (\ref{W10}) and $\sum_{n\in\mathbb{Z}}c(n)<+\infty$  that
\begin{eqnarray}
\label{W20}
\lim_{k\rightarrow +\infty}\sum_{n\in\mathbb{Z}}c(n)|u_k(n)|^{q-2}u_k(n)(u_k(n)-u(n))\ln{|u_k(n)|^r}=0.
\end{eqnarray}
Then using H\"{o}lder inequality
$$
d_1d_2+d_3d_4\leq(d_1^p+d_3^p)^{\frac{1}{p}}(d_2^{p^*}+d_4^{p^*})^{\frac{1}{{p^*}}},
$$
where $d_1, d_2, d_3, d_4$ are nonnegative constants and $p^*=\frac{p}{p-1}, p>1$, in virtue of (\ref{eq10}) and (\ref{eq6}), we have
\begin{small}
\begin{eqnarray}
\label{W16}
&    & \langle I'(u_k),u_k-u\rangle-\langle I'(u),u_k-u\rangle\nonumber\\
& =  & \sum_{n\in\mathbb{Z}}a(n)|\Delta u_k(n)|^{p-2}\Delta u_k(n)\Delta (u_k(n)-u(n))
          +\sum_{n\in\mathbb{Z}}b(n)|u_k(n)|^{p-2}u_k(n)(u_k(n)-u(n))\nonumber\\
&    & -\sum_{n\in\mathbb{Z}}c(n)|u_k(n)|^{q-2}u_k(n)(u_k(n)-u(n))\ln{|u_k(n)|^r}
          -\sum_{n\in\mathbb{Z}}a(n)|\Delta u(n)|^{p-2}\Delta u(n)\Delta (u_k(n)-u(n))\nonumber\\
&    & -\sum_{n\in\mathbb{Z}}b(n)|u(n)|^{p-2}u(n)(u_k(n)-u(n))
          +\sum_{n\in\mathbb{Z}}c(n)|u(n)|^{q-2}u(n)(u_k(n)-u(n))\ln{|u(n)|^r}\nonumber\\
& =  & \|u_k\|^p+\|u\|^p-\sum_{n\in\mathbb{Z}}a(n)|\Delta u_k(n)|^{p-2}\Delta u_k(n)\Delta u(n)
          -\sum_{n\in\mathbb{Z}}b(n)|u_k(n)|^{p-2}u_k(n)u(n)\nonumber\\
&    & -\sum_{n\in\mathbb{Z}}c(n)|u_k(n)|^{q-2}u_k(n)(u_k(n)-u(n))\ln{|u_k(n)|^r}
          -\sum_{n\in\mathbb{Z}}a(n)|\Delta u(n)|^{p-2}\Delta u(n)\Delta u_k(n)\nonumber\\
&    & -\sum_{n\in\mathbb{Z}}b(n)|u(n)|^{p-2}u(n)u_k(n)
          +\sum_{n\in\mathbb{Z}}c(n)|u(n)|^{q-2}u(n)(u_k(n)-u(n))\ln{|u(n)|^r}\nonumber\\
&\geq& \|u_k\|^p+\|u\|^p-\sum_{n\in\mathbb{Z}}a(n)|\Delta u_k(n)|^{p-1}|\Delta u(n)|
          -\sum_{n\in\mathbb{Z}}b(n)|u_k(n)|^{p-1}|u(n)|\nonumber\\
&    & -\sum_{n\in\mathbb{Z}}c(n)|u_k(n)|^{q-2}u_k(n)(u_k(n)-u(n))\ln{|u_k(n)|^r}
          -\sum_{n\in\mathbb{Z}}a(n)|\Delta u(n)|^{p-1}|\Delta u_k(n)|\nonumber\\
&    & -\sum_{n\in\mathbb{Z}}b(n)|u(n)|^{p-1}|u_k(n)|
          +\sum_{n\in\mathbb{Z}}c(n)|u(n)|^{q-2}u(n)(u_k(n)-u(n))\ln{|u(n)|^r}\nonumber\\
&\geq& \|u_k\|^p+\|u\|^p-\left(\sum_{n\in\mathbb{Z}}\big(a(n)^{\frac{1}{p}}|\Delta u(n)|\big)^p\right)^{\frac{1}{p}}
\left(\sum_{n\in\mathbb{Z}}\big(a(n)^{1-\frac{1}{p}}|\Delta u_k(n)|^{p-1}\big)^{\frac{p}{p-1}}\right)^{\frac{p-1}{p}}\nonumber\\
&    & -\left(\sum_{n\in\mathbb{Z}}\big(b(n)^{\frac{1}{p}}|u(n)|\big)^p\right)^{\frac{1}{p}}
\left(\sum_{n\in\mathbb{Z}}\big(b(n)^{1-\frac{1}{p}}|u_k(n)|^{p-1}\big)^{\frac{p}{p-1}}\right)^{\frac{p-1}{p}}\nonumber\\
&    & -\left(\sum_{n\in\mathbb{Z}}\big(a(n)^{\frac{1}{p}}|\Delta u_k(n)|\big)^p\right)^{\frac{1}{p}}
\left(\sum_{n\in\mathbb{Z}}\big(a(n)^{1-\frac{1}{p}}|\Delta u(n)|^{p-1}\big)^{\frac{p}{p-1}}\right)^{\frac{p-1}{p}}\nonumber\\
&    & -\left(\sum_{n\in\mathbb{Z}}\big(b(n)^{\frac{1}{p}}|u_k(n)|\big)^p\right)^{\frac{1}{p}}
\left(\sum_{n\in\mathbb{Z}}\big(b(n)^{1-\frac{1}{p}}|u(n)|^{p-1}\big)^{\frac{p}{p-1}}\right)^{\frac{p-1}{p}}\nonumber\\
&    & -\sum_{n\in\mathbb{Z}}c(n)|u_k(n)|^{q-2}u_k(n)(u_k(n)-u(n))\ln{|u_k(n)|^r}
          +\sum_{n\in\mathbb{Z}}c(n)|u(n)|^{q-2}u(n)(u_k(n)-u(n))\ln{|u(n)|^r}\nonumber\\
& =  & \|u_k\|^p+\|u\|^p-\left(\sum_{n\in\mathbb{Z}}a(n)|\Delta u(n)|^p\right)^{\frac{1}{p}}
\left(\sum_{n\in\mathbb{Z}}a(n)|\Delta u_k(n)|^{p}\right)^{\frac{p-1}{p}}\nonumber\\
&    & -\left(\sum_{n\in\mathbb{Z}}b(n)|u(n)|^p\right)^{\frac{1}{p}}
\left(\sum_{n\in\mathbb{Z}}b(n)|u_k(n)|^{p}\right)^{\frac{p-1}{p}}
       -\left(\sum_{n\in\mathbb{Z}}a(n)|\Delta u_k(n)|^p\right)^{\frac{1}{p}}
\left(\sum_{n\in\mathbb{Z}}a(n)|\Delta u(n)|^{p}\right)^{\frac{p-1}{p}}\nonumber\\
&    & -\left(\sum_{n\in\mathbb{Z}}b(n)|u_k(n)|^p\right)^{\frac{1}{p}}
\left(\sum_{n\in\mathbb{Z}}b(n)|u(n)|^{p}\right)^{\frac{p-1}{p}}
       -\sum_{n\in\mathbb{Z}}c(n)|u_k(n)|^{q-2}u_k(n)(u_k(n)-u(n))\ln{|u_k(n)|^r}\nonumber\\
&    &  +\sum_{n\in\mathbb{Z}}c(n)|u(n)|^{q-2}u(n)(u_k(n)-u(n))\ln{|u(n)|^r}\nonumber\\
&\geq& \|u_k\|^p+\|u\|^p-\left(\sum_{n\in\mathbb{Z}}[a(n)|\Delta u(n)|^p+b(n)|u(n)|^p]\right)^{\frac{1}{p}}
\left(\sum_{n\in\mathbb{Z}}[a(n)|\Delta u_k(n)|^{p}+b(n)|u_k(n)|^p]\right)^{\frac{p-1}{p}}\nonumber\\
&    & -\left(\sum_{n\in\mathbb{Z}}[a(n)|\Delta u_k(n)|^p+b(n)|u_k(n)|]\big)^p\right)^{\frac{1}{p}}
\left(\sum_{n\in\mathbb{Z}}[a(n)|\Delta u(n)|^{p}+b(n)|u(n)|]\big)^p\right)^{\frac{p-1}{p}}\nonumber\\
&    &  -\sum_{n\in\mathbb{Z}}c(n)|u_k(n)|^{q-2}u_k(n)(u_k(n)-u(n))\ln{|u_k(n)|^r}
           +\sum_{n\in\mathbb{Z}}c(n)|u(n)|^{q-2}u(n)(u_k(n)-u(n))\ln{|u(n)|^r}\nonumber\\
& =  & \|u_k\|^p+\|u\|^p-\|u\|\|u_k\|^{p-1}-\|u_k\|\|u\|^{p-1}
           -\sum_{n\in\mathbb{Z}}c(n)|u_k(n)|^{q-2}u_k(n)(u_k(n)-u(n))\ln{|u_k(n)|^r}\nonumber\\
&    &  +\sum_{n\in\mathbb{Z}}c(n)|u(n)|^{q-2}u(n)(u_k(n)-u(n))\ln{|u(n)|^r}\nonumber\\
& =  & \Big(\|u_k\|^{p-1}-\|u\|^{p-1}\Big)\Big(\|u_k\|-\|u\|\Big)
           -\sum_{n\in\mathbb{Z}}c(n)|u_k(n)|^{q-2}u_k(n)(u_k(n)-u(n))\ln{|u_k(n)|^r}\nonumber\\
&    &  +\sum_{n\in\mathbb{Z}}c(n)|u(n)|^{q-2}u(n)(u_k(n)-u(n))\ln{|u(n)|^r}.
\end{eqnarray}
\end{small}According to (\ref{W11}), (\ref{W15}), (\ref{W13}), (\ref{W20}) and (\ref{W16}),  we have $\|u_k\|\to \|u\|$ as $k\to +\infty$. According to the uniform convexity of $\mathcal{D}$ (similar to the argument of the Appendix A.1 in \cite{Yang 2023}), the fact that $u_k\rightharpoonup u$ in $\mathcal{D}$ and the Kadec-Klee property, we can obtain that $u_k\to u$ in $\mathcal{D}$. Thus, $I$ satisfies the Cerami condition.\qed
\par
Next, we prove that the functional $I$ defined by (\ref{eq7}) has a mountain pass geometry.
\vskip2mm
\noindent
{\bf Lemma 4.3.} {\it(i) There are two positive constants $\rho$ and $\delta'$ such that $I(u)\geq \delta'$  for all $u\in \mathcal{D}$ with $\|u\|=\rho$,\\
(ii) There is $\varphi_j\in \mathcal{D}\setminus \{0\}$ such that $I(t\varphi_j)\rightarrow -\infty$ as $t\rightarrow+\infty$.}\\
{\bf Proof.} For (i), it follows from $(\ref{eq7})$, $(\ref{eq9})$ and Lemma 2.1 that there exists $\varepsilon_7\in(0,\frac{qb_0}{pc_0})$ and $C_{\varepsilon_7}>0$ such that
\begin{eqnarray*}
       I(u)
&\geq& \frac{1}{p}\|u\|^p-\frac{1}{q}\sum_{n\in\mathbb{Z}}c(n)|u(n)|^q\ln {|u(n)|^r}\\
&\geq& \frac{1}{p}\|u\|^p-\frac{1}{q}c_0\varepsilon_7\|u\|_{l^p}^p
        -\frac{1}{q}c_0C_{\varepsilon_7}\|u\|^\zeta_{l^\zeta} \\
&\geq& \left(\frac{1}{p}-\frac{1}{q}c_0\varepsilon_7b_0^{-1}\right)\|u\|^p
        -\frac{1}{q}c_0C_{\varepsilon_7}b_0^{-\frac{\zeta}{p}}\|u\|^\zeta.
\end{eqnarray*}
Choose $\rho>0$ sufficiently small. There appears a constants $\beta_0=\frac{\rho^p q-\rho^p pc_0\left(\varepsilon_7b_0^{-1}+C_{\varepsilon_7}b_0^{-\frac{\zeta}{p}}\right)}{pq}>0$ such that $I(u)\geq\beta_0$ for all $u\in \mathcal{D}$ with $\|u\|=\rho$.
\par
By the definition of $c_*$, for any $j>0$,  we can choose a $\varphi_j\in \mathcal{N}\subset\mathcal{D}\setminus \{0\}$ such that
\begin{eqnarray}\label{pp1}
c_*\le I(\varphi_j)\le c_*+\frac{1}{j}.
\end{eqnarray}
Then for any  $t>0$, there holds
\begin{eqnarray*}
        I(t\varphi_j)
& =  &\frac{t^p}{p}\|\varphi_j\|^p+\frac{rt^q}{q^2}\sum_{n\in\mathbb{Z}}c(n)|\varphi_j(n)|^q
        -\frac{t^q\ln {t^r}}{q}\sum_{n\in\mathbb{Z}}c(n)|\varphi_j(n)|^q
        -\frac{t^q}{q}\sum_{n\in\mathbb{Z}}c(n)|\varphi_j(n)|^q\ln {|\varphi_j(n)|^r}\\
&\leq &\frac{t^p}{p}\|\varphi_j\|^p+\frac{rt^q}{q^2}\sum_{n\in\mathbb{Z}}c(n)|\varphi_j(n)|^q
        -\frac{t^q\ln {t^r}}{q}\sum_{n\in\mathbb{Z}}c(n)|\varphi_j(n)|^q
          +\frac{t^q}{q}\sum_{\substack{n\in\mathbb{Z}}}c_0(\varepsilon|\varphi_j(n)|^p+C_\varepsilon|\varphi_j(n)|^\zeta)
\end{eqnarray*}
which implies that $I(t\varphi_j)\rightarrow -\infty$ as $t\rightarrow +\infty$, and so (ii) holds.\qed
\vskip2mm
\par
\vskip2mm
\noindent
{\bf Proof of Theorem 1.2.}  Lemma 4.1 and Lemma 4.3 imply that $I$ has a Cerami sequence $\{u_{kj}\}$ at the level $d_j$, that is,
\begin{eqnarray*}
\label{WW1}
I(u_{kj})\rightarrow d_j\;\;\text{ and}\;\; \|I'(u_{kj})\|(1+\|u_{kj}\|)\rightarrow 0, \text{ as }kj\to+\infty.
\end{eqnarray*}
In virtue of Remark 4.1 and the definition of $d_j$, it is easy to see that $d_j\in [\beta_0, \max\limits_{0\leq t\leq 1} I(t\varphi_j)]$. Furthermore, noting that $\varphi_j\in \mathcal{N}$, according to Corollary 2.5, we obtain that $I(\varphi_j)=\max_{t\ge 0} I(t\varphi_j)$, and so  $d_j\in [\beta_0, I(\varphi_j)]$ which together with (\ref{pp1}) implies that $d_j\in [\beta_0, c_*+\frac{1}{j}]$. Thus, we can choose a subsequence $\{u_{k_j,j}\}$, denoted by $\{u_j\}$, such that
\begin{eqnarray}
\label{pp3}
I(u_{j})\rightarrow d_*\;\;\text{ and}\;\; \|I'(u_{j})\|(1+\|u_{j}\|)\rightarrow 0, \text{ as }j\to+\infty,
\end{eqnarray}
for some $d_*\in [\frac{\beta_0}{2}, c_*].$ (\ref{pp3}) and Lemma 4.2 imply that $\{u_{j}\}$ has a convergent subsequence, still denoted by $\{u_{j}\}$, such that $u_j\to \bar{u}$ as $j\to+\infty$.
By the continuity of $I$ and $I'$, we obtain that  $I(\bar{u})=d_*$ and  $I'(\bar{u})=0$, which together with the fact $d_*\ge \frac{\beta_0}{2}>0$ implies that $\bar{u}\in\mathcal{N}$ is a nontrivial solution of $(\ref{eq1})$ and obviously, $I(\bar{u})\geq c_*=\inf_{u\in\mathcal{N}}I(u)$.  Moreover, according to (\ref{pp3}), $(\ref{eq7})$, $(\ref{eq10})$ and the weak lower semi-continuity of norm, there exists
\begin{eqnarray*}
       c_*
& \ge & d_*=\lim_{j\to+\infty}\left[I(u_j)-\frac{1}{p}\langle I'(u_j),u_j\rangle\right]
              =\lim_{j\to +\infty}\left[\left(\frac{1}{p}-\frac{1}{q}\right)\|u_j\|^p
               +\frac{r}{q^2}\sum_{n\in\mathbb{Z}}c(n)|u_j(n)|^q\right]\\
&\geq& \left(\frac{1}{p}-\frac{1}{q}\right)\|\bar{u}\|^p+\frac{r}{q^2}\sum_{n\in\mathbb{Z}}c(n)|\bar{u}(n)|^q
              =I(\bar{u})-\frac{1}{p}\langle I'(\bar{u}),\bar{u}\rangle= I(\bar{u}),
\end{eqnarray*}
which implies that $I(\bar{u})\leq c_*$. Thus $I(\bar{u})=c_*=\inf_{u\in\mathcal{N}} I(u)>0$.
\par
Finally, we prove that $m_*\ge 2c_*$. In fact, it follows from Corollary 2.4, Proposition 2.1, Lemma 2.6 and Lemma 2.8 that there are two positive constants $s'$ and $t'$ such that $s'u_0^+\in\mathcal{N}$ and $t'u_0^-\in\mathcal{N}$, then one has
\begin{eqnarray*}
        m_*=I(u_0)
&=&\max_{s,t\geq0}I(su_0^++tu_0^-)\\
&\geq&  \max_{s,t\geq0}\left[I(su^+_0)+I(tu_0^-)\right]\\
              &=&\max_{s\geq0}I(su^+_0)+\max_{t\geq0}I(tu^-_0)\\
&\geq&  I(s'u^+_0)+I(t'u^-_0)\\
&\geq&  I(\bar{u})+I(\bar{u})=2c_*.
\end{eqnarray*}
The proof of Theorem 1.2 is completed.\qed

\vskip2mm
{\section{Appendix}}
\noindent{\bf Appendix 1.} {\it There exists $u\in E$ such that $\sum\limits_{n\in\mathbb{Z}}c(n)|u(n)|^q\ln{|u(n)|^r}=-\infty$.}\\
{\bf Proof.} Set
\begin{eqnarray*}
   u(n)
=
   \begin{cases}
   \frac{1}{|n|\ln{|n|}},&|n|\geq p+2,\\
   0,&|n|\leq p+1,
   \end{cases}
\end{eqnarray*}
where $|n|$ represents the absolute value of $n$. Let
\begin{eqnarray*}
   a(n),b(n)
=
   \begin{cases}
   |n|^{p-1},&|n|\geq p-1,\\
   1,&|n|\leq p-2.
   \end{cases}
\;\;\text{and}\;\;
   c(n)
=
   \begin{cases}
   |n|^{q-1},&|n|\geq q-1,\\
   1,&|n|\leq q-2.
   \end{cases}
\end{eqnarray*}
Using integral discriminant method and convergence and divergence of improper integral, we can obtain the following results
\begin{eqnarray*}
   \sum_{n=2}^{+\infty}\frac{1}{n(\ln n)^\theta}
=
   \begin{cases}
   < +\infty,\;\text{that is, the series convergence,}& \theta>1,\\
   = +\infty,\;\text{that is, the series divergence,}& 0\leq\theta\leq1.
   \end{cases}
\end{eqnarray*}
According to the definitions of $\Delta$, $u(n)$ and $a(n)$, using $C_p$ inequality, we have
\begin{small}
\begin{eqnarray*}
       \sum_{n\in\mathbb{Z}}a(n)|\Delta u(n)|^p
& =  & |p+1|^{p-1}\left(\frac{1}{|p+2|\ln{|p+2|}}\right)^p
         +\sum_{|n|\geq p+2}|n|^{p-1}\left|\frac{1}{|n+1|\ln{|n+1|}}-\frac{1}{|n|\ln{|n|}}\right|^p\\
&\leq& \frac{1}{|p+2|(\ln{|p+2|})^p}
         +\sum_{|n|\geq p+2}|n|^{p-1}\left(\frac{1}{|n+1|\ln{|n+1|}}+\frac{1}{|n|\ln{|n|}}\right)^p\\
&\leq& \frac{1}{|p+2|(\ln{|p+2|})^p}+\sum_{|n|\geq p+2}|n|^{p-1}
                2^{p-1}\left[\left(\frac{1}{|n+1|\ln{|n+1|}}\right)^p+\left(\frac{1}{|n|\ln{|n|}}\right)^p\right]\\
&\leq& \frac{1}{|p+2|(\ln{|p+2|})^p}+2^{p-1}\sum_{|n|\geq p+2}\left[\frac{1}{|n+1|(\ln{|n+1|})^p}
            +\frac{1}{|n|(\ln{|n|})^p}\right]\\
& <  & +\infty.
\end{eqnarray*}
\end{small}
Then according to the definition of $u(n)$ and $b(n)$, there holds
$$
  \sum_{n\in\mathbb{Z}}b(n)|u(n)|^p
= \sum_{|n|\geq p+2}\frac{|n|^{p-1}}{(|n|\ln {|n|})^p}
= \sum_{|n|\geq p+2}\frac{1}{|n|(\ln{|n|})^p}
< +\infty.
$$
Therefore, it is easy to know that $u\in E$.
\par
Now, we prove that $\sum_{n\in\mathbb{Z}}c(n)|u(n)|^q\ln{|u(n)|^r}=-\infty$ if $1<q\leq 2$. Note that if $|n|\ge p+2$, then
$$
  |u(n)|^q\ln{|u(n)|^r}
= \frac{1}{(|n|\ln{|n|})^q}\times \ln{\frac{1}{(|n|\ln{|n|})^r}}
= -\left(\frac{r}{|n|^q(\ln{|n|})^{q-1}}+\frac{r\ln(\ln{|n|})}{|n|^q(\ln{|n|})^q}\right).
$$
Thus, there holds
\begin{eqnarray*}
       \sum_{n\in\mathbb{Z}}c(n)|u(n)|^q\ln{|u(n)|^r}
& = & \sum_{|n|\geq p+2}|n|^{q-1}\times
         \left[-\left(\frac{r}{|n|^q(\ln{|n|})^{q-1}}+\frac{r\ln(\ln{|n|})}{|n|^q(\ln{|n|})^q}\right)\right]\\
& = & -\left(\sum_{|n|\geq p+2}\frac{r}{|n|(\ln{|n|})^{q-1}}
              +\sum_{|n|\geq p+2}\frac{r\ln(\ln{|n|})}{|n|(\ln{|n|})^q}\right)\\
& := & -(I'+II').
\end{eqnarray*}
Because $|n|\geq p+2$,  we have  $\ln \ln |n|>0$ and then $II'>0$. On the other hand, it is easy to see $I'=\infty$ since $1<q\leq 2$. Thus we complete the proof.
\qed
\vskip2mm
\noindent
{\bf Appendix 2.} {\it For all $s,t\geq0$ and $0\leq i,j\leq\frac{p}{2}$, there exist\\
$\frac{2s^pC_{\frac{p}{2}-1}^iC_i^j+s^pC_{\frac{p}{2}-1}^{i-1}C_{i-1}^j+2t^pC_{\frac{p}{2}-1}^{i-1}C_{i-1}^{j-1}
    +t^pC_{\frac{p}{2}-1}^{i-1}C_{i-1}^j-2s^{p-(i+j)}t^{i+j}C_{\frac{p}{2}}^i C_i^j}{2p}\geq0$,
$\frac{s^p C_{\frac{p}{2}-1}^{i}+t^pC_{\frac{p}{2}-1}^{i-1}-s^{p-2i}t^{2i}C_{\frac{p}{2}}^{i}}{p}\geq0$,
$\frac{2s^pC_{\frac{p}{2}-1}^{i}+s^pC_{\frac{p}{2}-1}^{i-1}+t^pC_{\frac{p}{2}-1}^{i-1}-2s^{p-i}t^iC_{\frac{p}{2}}^{i}}{2p}\geq0$
\;\;and\;\;\;$\frac{s^pC_{\frac{p}{2}-1}^{j}+2t^pC_{\frac{p}{2}-1}^{j-1}+t^pC_{\frac{p}{2}-1}^{j}
-2s^{\frac{p}{2}-j}t^{\frac{p}{2}+j}C_{\frac{p}{2}}^{j}}{2p}\geq0$.}\\
{\bf Proof.} Using combination formula $C_{p-1}^j=\frac{p-j}{p}C_p^j$, $C_{p-1}^{j-1}=\frac{j}{p}C_p^j$ and deformation of Young inequality $x^\lambda y^{1-\lambda}\leq\lambda x+(1-\lambda)y,\;(0<\lambda<1)$, there exists
\begin{eqnarray*}
&    & \frac{2s^pC_{\frac{p}{2}-1}^iC_i^j+s^pC_{\frac{p}{2}-1}^{i-1}C_{i-1}^j+2t^pC_{\frac{p}{2}-1}^{i-1}C_{i-1}^{j-1}
           +t^pC_{\frac{p}{2}-1}^{i-1}C_{i-1}^j-2s^{p-(i+j)}t^{i+j}C_{\frac{p}{2}}^i C_i^j}{2p}\\
&  = & \frac{2s^p(\frac{\frac{p}{2}-i}{\frac{p}{2}}C_{\frac{p}{2}}^i)C_i^j+s^p(\frac{i}{\frac{p}{2}}C_{\frac{p}{2}}^{i})
          (\frac{i-j}{i}C_{i}^j)+2t^p(\frac{i}{\frac{p}{2}}C_{\frac{p}{2}}^{i})(\frac{j}{i}C_{i}^{j})+t^p(\frac{i}{\frac{p}{2}}
          C_{\frac{p}{2}}^{i})(\frac{i-j}{i}C_{i}^j)-2s^{p-(i+j)}t^{i+j}C_{\frac{p}{2}}^i C_i^j}{2p}\\
&  = & \frac{C_{\frac{p}{2}}^iC_i^j\left[2(p-2i)s^p+2(i-j)s^p+4jt^p+2(i-j)t^p-2ps^{p-(i+j)}t^{i+j}\right]}{2p^2}\\
&  = & \frac{C_{\frac{p}{2}}^iC_i^j}{p}\left[\frac{p-(i+j)}{p}s^p+\frac{i+j}{p}t^p-s^{p-(i+j)}t^{i+j}\right]\\
&\geq& \frac{C_{\frac{p}{2}}^iC_i^j}{p}\left[(s^p)^{\frac{p-(i+j)}{p}}(t^p)^{\frac{i+j}{p}}-s^{p-(i+j)}t^{i+j}\right]=0.
\end{eqnarray*}
Similarly, we can obtain that $\frac{2s^pC_{\frac{p}{2}-1}^{i}+s^pC_{\frac{p}{2}-1}^{i-1}+t^pC_{\frac{p}{2}-1}^{i-1}-2s^{p-i}t^iC_{\frac{p}{2}}^{i}}{2p}\geq0$, $\frac{s^p C_{\frac{p}{2}-1}^{i}+t^pC_{\frac{p}{2}-1}^{i-1}-s^{p-2i}t^{2i}C_{\frac{p}{2}}^{i}}{p}\geq0$ and
$\frac{s^pC_{\frac{p}{2}-1}^{j}+2t^pC_{\frac{p}{2}-1}^{j-1}+t^pC_{\frac{p}{2}-1}^{j}
-2s^{\frac{p}{2}-j}t^{\frac{p}{2}+j}C_{\frac{p}{2}}^{j}}{2p}\geq0$.
Thus, we can obtain that the conclusions hold.
\qed
\vskip2mm
\noindent
{\bf Appendix 3.} {\it For all $\frac{s_2}{s_1},\frac{t_2}{t_1}\geq0$ and $0\leq i,j\leq\frac{p}{2}$, there exists\\
$$
\frac{2(\frac{s_2}{s_1})^pC_{\frac{p}{2}-1}^iC_i^j+(\frac{s_2}{s_1})^pC_{\frac{p}{2}-1}^{i-1}C_{i-1}^j
+2(\frac{t_2}{t_1})^pC_{\frac{p}{2}-1}^{i-1}C_{i-1}^{j-1}+(\frac{t_2}{t_1})^pC_{\frac{p}{2}-1}^{i-1}C_{i-1}^j
-2(\frac{s_2}{s_1})^{p-(i+j)}(\frac{t_2}{t_1})^{i+j}C_{\frac{p}{2}}^iC_i^j}{2p}\geq0.
$$}
{\bf Proof.} Using the method of Appendix 2, we have
\begin{eqnarray*}
&    & \frac{2(\frac{s_2}{s_1})^pC_{\frac{p}{2}-1}^iC_i^j+(\frac{s_2}{s_1})^pC_{\frac{p}{2}-1}^{i-1}C_{i-1}^j
+2(\frac{t_2}{t_1})^pC_{\frac{p}{2}-1}^{i-1}C_{i-1}^{j-1}+(\frac{t_2}{t_1})^pC_{\frac{p}{2}-1}^{i-1}C_{i-1}^j
-2(\frac{s_2}{s_1})^{p-(i+j)}(\frac{t_2}{t_1})^{i+j}C_{\frac{p}{2}}^iC_i^j}{2p}\\
&  = & \frac{C_{\frac{p}{2}}^iC_i^j\left[2(p-2i)(\frac{s_2}{s_1})^p+2(i-j)(\frac{s_2}{s_1})^p
           +4j(\frac{t_2}{t_1})^p+2(i-j)(\frac{t_2}{t_1})^p-2p(\frac{s_2}{s_1})^{p-(i+j)}(\frac{t_2}{t_1})^{i+j}\right]}{2p^2}\\
&  = & \frac{C_{\frac{p}{2}}^iC_i^j}{p}\left[\frac{p-(i+j)}{p}\left(\frac{s_2}{s_1}\right)^p
           +\frac{i+j}{p}\left(\frac{t_2}{t_1}\right)^p-\left(\frac{s_2}{s_1}\right)^{p-(i+j)}\left(\frac{t_2}{t_1}\right)^{i+j}\right]\\
&\geq& \frac{C_{\frac{p}{2}}^iC_i^j}{p}\left[\left(\left(\frac{s_2}{s_1}\right)^p\right)^{\frac{p-(i+j)}{p}}
         \left(\left(\frac{t_2}{t_1}\right)^p\right)^{\frac{i+j}{p}}-\left(\frac{s_2}{s_1}\right)^{p-(i+j)}\left(\frac{t_2}{t_1}\right)^{i+j}\right]=0.
\end{eqnarray*}
Thus we complete the proof.
\qed
\vskip3mm
\noindent{\bf Funding information}
\par
\noindent
This project is  supported by Yunnan Fundamental Research Projects (grant No: 202301AT070465) and supported by Yunnan Ten Thousand Talents Plan Young \& Elite Talents Project.

\vskip3mm
 \noindent
\noindent{\bf Author contribution}
\par
\noindent
Ou and Zhang contributed equally to this work. All authors reviewed the manuscript.

\vskip3mm
 \noindent
\noindent{\bf Conflict of interest}

\noindent
The authors state no conflict of interest.

\vskip3mm
 \noindent
\noindent{\bf Data availability statement }

\noindent
Not available.

\vskip2mm
\renewcommand\refname{References}
{}
 \end{CJK}
\end{document}